\documentclass[a4paper,11pt,reqno,english]{smfart}
\usepackage{hyperref} 
\usepackage{psfrag}
\usepackage{latexsym}
\usepackage{array}
\usepackage{amssymb}
\usepackage{smfthm}
\theoremstyle{plain}

\usepackage{amsmath,amssymb,amsfonts,graphicx,psfrag,subfig}
\usepackage[latin1]{inputenc}
\newtheorem{thm}{Theorem}[section]
\newtheorem{lem}[thm]{Lemma}
\newtheorem{df}[thm]{Definition}

\newtheorem{rem}[thm]{Remark}

\newtheorem{cor}[thm]{Corollary}

\def\be#1 {\begin{equation} \label{#1}}
\newcommand{\ee}{\end{equation}}
\def\dem {\noindent {\bf Proof: }}

\newcommand{\mb}{\medskip\noindent}
\newcommand{\gb}{\bigskip\noindent}
\newcommand{\R}{\mathbb R}

\newcommand{\B}{\mathbb{B}}

\newcommand{\F}{\mathcal F}

\newcommand{\E}{\mathbb E}

\newcommand{\espace}{(\Omega,\F,\PP)}
\newcommand{\I}{\mathcal I}

\newcommand{\PP}{\mathbb P}
\newcommand{\PPP}{\mathrm P}
\def \Qc {\tilde{Q}}

\def \qqq {{\bf q}}
\def \qq {\mathrm{q}}

\def \UU {{\bf U}}
\def \BB {{\mathcal{B}}}

\def \NN {\mathrm{N}}

\def \virg {\, , \,\,}
\def \dsp {\displaystyle}
\def \vsp {\vspace{6pt}}

\def\sqw{\hbox{\rlap{\leavevmode\raise.3ex\hbox{$\sqcap$}}$%
\sqcup$}}
\def\findem{\ifmmode\sqw\else{\ifhmode\unskip\fi\nobreak\hfil
\penalty50\hskip1em\null\nobreak\hfil\sqw
\parfillskip=0pt\finalhyphendemerits=0\endgraf}\fi}

\begin{document}

\title[Stochastic sweeping process]{Stochastic perturbation of sweeping process and a convergence result for an associated numerical scheme}
\date{January 18, 2010}

\author{Fr\'ed\'eric Bernicot}
\address{CNRS - Universit\'e Lille 1 \\ Laboratoire de math\'ematiques Paul Painlev\'e \\ 59655 Villeneuve d'Ascq Cedex, France}
\email{frederic.bernicot@math.univ-lille1.fr} \urladdr{http://math.univ-lille1.fr/$\sim$bernicot/}
\author{Juliette Venel}
\address{LAMAV \\ Universit\'e de Valenciennes et du Hainaut-Cambr\'esis\\Mont Houy 59313 Valenciennes Cedex 9 }
\email{juliette.venel@univ-valenciennes.fr} \urladdr{http://www.math.u-psud.fr/$\sim$venel/}

\frontmatter 
\begin{abstract} Here we present well-posedness results for first order stochastic differential inclusions, more precisely for sweeping process with a stochastic perturbation. These results are provided in combining both deterministic sweeping process theory (recently developed in \cite{Thibrelax} and \cite{Thibbv}) and methods concerning the reflection of a Brownian motion (\cite{LS} and \cite{Saisho}). In addition, we prove convergence results for a Euler scheme, discretizing theses stochastic differential inclusions.
\end{abstract}

\begin{altabstract} Nous d\'emontrons dans ce travail le caract\`ere ``bien-pos\'e'' d'inclusions diff\'erentielles stochastiques du premier ordre, plus pr\'ecis\'ement de processus de rafle avec une perturbation stochastique. Ces r\'esultats sont issus de l'association de la th\'eorie des processus de rafle d\'eterministes (r\'ecemment d\'evelopp\'ee \cite{Thibrelax} et \cite{Thibbv}) et de m\'ethodes concernant la r\'eflexion d'un mouvement brownien (\cite{LS} et \cite{Saisho}). De plus, nous prouvons un r\'esultat de convergence pour un sch\'ema d'Euler, discr\'etisant ces inclusions diff\'erentielles stochastiques.
\end{altabstract}

\subjclass{34A60 ; 65L20 ; 60H10}
\keywords{Sweeping process ; differential inclusions ; stochastic differential equations ; Euler scheme}
\altkeywords{ Processus de rafle ; inclusions diff\'erentielles ; \'equations diff\'erentielles stochastiques ; Sch\'ema d'Euler}

\maketitle
\mainmatter

\tableofcontents

\gb

\section{Introduction}

In this paper, we are interested in particular first order differential inclusions (namely sweeping process) with a stochastic perturbation. This work rests on the combining of two different theories: the first one about sweeping process and the second one about the reflection of a Brownian motion on a boundary. Let us first recall these two problems and related results.

\gb {\bf Sweeping process}

Let ${\mathcal B}$ be a Banach space, $\I$ be a bounded time-interval, $C:\I \rightrightarrows {\mathcal B}$ be a set-valued map with nonempty closed values, and let $f: \I \times {\mathcal B} \rightrightarrows {\mathcal B}$ be a perturbation. The associated sweeping process takes the form:
\begin{equation} 
\left\{
\begin{array}{l}
 \dsp \frac{du}{dt}(t) + \NN(C(t),u(t)) \ni f(t,u(t)) \vsp \\
 u(t)\in C(t) \vsp \\
 u(0)=u_0\ ,
\end{array}
\right. \label{sys1}
\end{equation}
with an initial data $u_0\in C(0)$ and where $\NN(C,x)$ denote the proximal normal cone of $C$ at any point $x$. This differential inclusion can be thought as follows: the point $u(t)$, submitted to the field $f(t,u(t))$, has to remain in the set $C(t)$.

\gb This type of evolution problem has been extensively studied. It has been introduced by
 J.J. Moreau in 70's  (see \cite{Moreausweep}) with convex sets $C(t)$ of a Hilbert space and with no perturbation ($f\equiv 0$). To solve this problem, J.J. Moreau brings a new important idea in proposing a \textit{catching-up algorithm}.

\mb Since then, some attempts have been made in the litterature to weaken the assumptions, for example to add a perturbation $f$, to weaken the convexity assumption of the sets, to obtain results in Banach spaces (and not only in Hilbert spaces). \\
The perturbed problem in finite dimension (${\mathcal B}=\R^d$) with convex sets $C(t)$ (or complements of convex sets) was firstly studied by C.~Castaing, T.X. D\'uc H\={a} and M. Valadier in \cite{CXDV}. In this framework, they proved existence of solutions for (\ref{sys1}) with a convex compact valued perturbation $f$ and a Lipschitzean set-valued map $C$.  Then in \cite{Castaingper}, C.~Castaing and M.D.P.~Monteiro Marques considered similar problems in assuming the upper semicontinuity of $f$ and a ``linear compact growth'':
\begin{equation}
 f(t,x)\subset \beta(t) (1+|x|)\overline{B(0,1)} \virg \forall (t,x) \in I \times \R^d.
\label{hypF}
\end{equation}
Moreover the set-valued map $C$ was supposed to be Hausdorff continuous and satisfying an ``interior ball condition'': 
\begin{equation} 
 \exists r>0 \virg  B(0,r) \subset C(t) \virg \forall t \in I.
\label{hypC}
\end{equation}

\gb  Later the main concept which weakens the convexity property of sets $C(t)$, is the notion of  ``uniform prox-regularity''. A set $C$ is said to be {\it uniformly prox-regular with constant $\eta$} or {\it $\eta$-prox-regular} if the projection onto $C$ is single-valued and continuous at any point distant at most $\eta$ from $C$. \\
The uniform prox-regularity assumption was made in numerous works dealing with sweeping process. The Hilbertian case without perturbation ($f\equiv 0$) was firstly treated by G.~Colombo, V.V.~Goncharov in \cite{Colombo}, by H.~Benabdellah in \cite{Benab} and later by L.~Thibault in \cite{Thibsweep} and by G. Colombo, M.D.P.~Monteiro Marques in \cite{Monteiro}. In \cite{Thibsweep}, the considered problem is 
\begin{equation}
 \left \{
\begin{array}{l} 
-du \in \NN(C(t), u(t)) \vsp \\
u(T_0)=u_0\ ,
 \end{array}
\right.
\label{eq:mesdiffssm}
\end{equation} 
where $du$ is the differential measure of $u$. The well-posedness of (\ref{eq:mesdiffssm}) is proved under the same assumptions as previously excepted (\ref{hypC}). \\
In an infinite dimensional Hilbert space ${\mathcal B}$ (${\mathcal B}=H$), the perturbed problem is studied by M.~Bounkhel, J.F.~Edmond and L.~Thibault in \cite{Thibnonconv, Thibsweep, Thibrelax, Thibbv} (see Theorem \ref{thm1}) and recently by the authors in a Banach space in \cite{BV}. For example in \cite{Thibbv}, the well-posedness of
\begin{equation}
 \left \{
\begin{array}{l} 
-du \in \NN(C(t), u(t)) + f(t,u(t)) dt \vsp \\
 u(0)=u_0\ ,
 \end{array}
\right.
\label{eq:mesdiffasm}
\end{equation} 
is proved with a set-valued map $C$ taking $\eta$-prox regular values (for some $\eta>0$) such that   
\begin{equation}
|d_{C(t)}(y) - d_{C(s)}(y)| \leq \mu(]s,t]) \virg \forall y\in H,\ \forall \  s, t \in I  \virg s\leq t
\label{varadon}
\end{equation}
where $\mu$ is a nonnegative measure satisfying
\begin{equation}
\sup_{s \in I} \mu(\{s\}) <\frac{\eta}{2}.
 \label{charge_singleton}
\end{equation}
All the proofs rest on the algorithm developed by J.J.~Moreau with additional arguments to deal with the prox-regularity assumption.

\gb We now want to recall results about the reflected Brownian motion in a set. Before that, we refer the reader to \cite{C1, C2}, where C. Castaing proved some results about existence of solutions for sweeping process with a convex set $C(t)$ which is stochastically perturbated. Here we want to add a stochastic perturbation in the differential inclusion without changing the deterministic time-evolution of $C$.

\gb {\bf Reflected Brownian motion}

We consider a closed set $C$ in $\R^d$ and we look for solving stochastic differential inclusions, describing the time-evolution of a Brownian motion (inside $C$) with a reflecting boundary $\partial C$. Let $(B_t)_{t>0}$ be a $\R^d$-valued Brownian motion, then the path is given by a stochastic process $X$ involving the following stochastic differential inclusion:
$$ 
\left\{
\begin{array}{l}
 \dsp dX_t + \NN(C,X_t) \ni dB_t \vsp \\
 X_t\in C \vsp \\
 X_0=u_0,
\end{array}
\right. 
$$
where $u_0$ is the starting point. One of the difficulty is to give a precise sense to ``$\NN(C,X_t)$''.

\mb The first well-posedness results have been investigated by many authors (see for example A.V.~Skorohod \cite{Skorohod1,Skorohod2}, N. Ikeda and S. Watanabe \cite{IW,W} and N. El Karoui \cite{ElKaroui} ...) for $C$ a half-plane (or a half-line). Then the problem was treated for a smooth set $C$ by D.W. Stroock and S.R.S. Varadhan \cite{SV} and later by A. Bensoussan and J.L. Lions \cite{BL}. H. Tanaka succeeded in getting around the smoothness assumption in the case of a convex set in \cite{T}. In \cite{LS}, P.L. Lions and A.S. Sznitman gave the first result for bounded uniformly prox-regular sets (without this terminology) in assuming an extra assumption of ``admissibility''. Few years later, Y. Saisho has extended the proof for unbounded prox-regular sets in \cite{Saisho}. They considered the associated deterministic Skorohod problem, which consists for a continuous function $h$ on a time-interval $\I$ in finding a pair of continuous functions $(x,k)$ defined on $\I$ with $k\in BV(\I)$ satisfying:
\be{eq2:intro} \forall t\in\I, \qquad x(t)+k(t)=h(t) \ee
and
\be{eq3:intro} |k|(t) = \int_0^t {\bf 1}_{x(s)\in\partial C} d|k|(s), \qquad k(t) = \int_0^t \xi_s d|k|(s) ,\ee
with $\xi_s\in \NN(C,x(s))$.\\
The equation (\ref{eq2:intro}) corresponds to an integral version of (\ref{sys1}) with a constant set $C(t):=C$ and $h(t)=\int_0^t f(s,u(s))ds$. Indeed (\ref{eq3:intro}) specifies that the support of the differential measure $dk$ is the set of the snapshots $t$ when $x$ reaches the boundary $\partial C$ and gives a precise sense to ``$dk\in \NN(C,x(t))$''. 
So the Skorohod problem can be thought as an integral version of a sweeping process.
Using this deterministic problem, the previously cited works deal with the following stochastic differential inclusion
$$ 
\left\{
\begin{array}{l}
 \dsp dx(t) + \NN(C,x(t)) \ni f(t,x(t))dt + \sigma(t,x(t)) dB_t \vsp \\
 u(t)\in C \vsp \\
 u(0)=u_0\ ,
\end{array}
\right. 
$$
where $f$ and $\sigma$ are some perturbations and $(B_t)_{t>0}$ is a real Brownian motion. This problem can also be seen as a stochastic perturbation of the ``sweeping process'' by a constant set $C$. We would like to finish this brief state of art by referring the reader to \cite{Saisho2} for a work about reflecting Brownian motion in a set corresponding to the complement of a collection of balls. We are specially interested by this example and we are looking for extending this result to more general situations with time-dependent constraints. 

\gb {\bf Associated numerical schemes}

Concerning the deterministic sweeping process, the existence results are obtained by the convergence of the so-called ``catching-up algorithm''. More precisely, in considering some subdivision $(J_k)_k$ of the time-interval, the set-valued map $C $ is approached by a piecewise constant multifunction taking value $C_k $ on $J_k=[t_k,t_{k+1}[$. The associated discretized solution $u $ is defined as follows: $$ \forall t \in J_{k+1} \virg   u(t)=u_{k+1} = \PPP_{C_{k+1}}\left[u_{k}+(t_{k+1}-t_k)f(t_k,u_k)\right],$$ 
with $u_0 $ fixed to the initial value.
Then it is proved that under the above assumptions, the sequence of discretized solutions is of Cauchy type and also converges to a function, which is the solution of the continuous problem (\ref{sys1}). 

\mb About the first order stochastic differential inclusions, we would like to present the work of E. Cepa \cite{Cepa} and F. Bernardin \cite{Bernardin}. They consider equations taking the form
$$ 
\left\{
\begin{array}{l}
 \dsp dx(t) + A(x(t))dt \ni f(t,x(t))dt + \sigma(t,x(t)) dB_t \vsp \\
 u(t)\in C \vsp \\
 u(0)=u_0\ ,
\end{array}
\right. 
$$
where $A$ is a time-independent maximal monotone operator. They define Euler numerical scheme (using the resolvents and Yosida's approximation of a monotone operator) and prove its convergence. In our case, (even for time-independent set) we are interested in uniformly prox-regular set $C$ and so, it is well-known that the associated proximal normal cone (as multivalued operator) could be not maximal monotone. It only satisfies a weaker property of hypomonotonicity.

\gb {\bf Framework for Stochastic perturbation of sweeping process and main results}

 We now come to our contribution. As explained before, the theories of sweeping process and reflected Brownian motion are based on the uniform prox-regularity property of the sets. In order to combine these two theories, we consider a set-valued map $C(\cdot)$ which is admissible and regular (in extending these defintions to set-valued maps, see Definitions \ref{def:adm} and \ref{def:reg}).

\mb We begin this work by describing some abstract results about uniformly prox-regular sets in order to understand our required assumption.
Then, we propose general results about a set-valued map $C$. Extending the assumptions of \cite{LS} to the time-dependent framework, we prove the two following results in Section \ref{sec:2}. 

\begin{thm} Let $\I$ be a bounded time interval and $C:\I\rightrightarrows \R^d$ admissible, regular set-valued map varying in an absolutely continuous way  (see Definitions \ref{def:adm}, \ref{def:reg} and \ref{def:abs}).\\
Let $f,\sigma:\I\times \R^d \rightarrow \R^d$ be bounded and Lipschitz with respect to the second variable.
Then the following problem
\begin{equation} 
\left\{
\begin{array}{l}
 \dsp dx(t) + \NN(C(t),x(t)) \ni f(t,x(t))dt + \sigma(t,x(t)) dB_t \vsp \\
 u(t)\in C(t) \vsp \\
 u(0)=u_0, 
\end{array}
\right. \label{eq:sto}
\end{equation}
(given by a $\R$-Brownian motion $(B_t)_{t\in \I}$ on a probability space $(\Omega,\F,{\mathbb P})$) is well-posed in $L^4(\Omega,L^\infty(\I))$. That means: for every initial data $u_0\in C(0)$, there exists one and only one process $(X_t)_{t\in \I}\in L^4(\Omega,L^\infty(\I))$ solution of (\ref{eqsto}) (in the sense of pathwise uniqueness). 
 \end{thm}

\noindent Moreover we will obtain stability results, which permit to understand the behavior of the stochastic solutions when the perturbation $\sigma$ tends to $0$ (see Theorem \ref{thm:conv1}).

\mb If the regularity of $C(\cdot)$ is not assumed, well-posedness results still hold (but not necessary in the space $L^4(\Omega,L^\infty(\I))$):

\begin{thm} Let $\I$ be a bounded time interval and $C:\I\rightrightarrows \R^d$ admissible  set-valued map varying in an absolutely continuous way  (see Definitions \ref{def:adm} and \ref{def:abs}).\\
Under the same assumptions on $f$ and $\sigma$, the problem (\ref{eq:sto}) has a unique solution (in the sense of pathwise uniqueness).
\end{thm}

\mb In Section \ref{sec:schema}, we prove convergence results for a Euler scheme, discretizing the stochastic differential inclusion (\ref{eq:sto}).

\mb In Section \ref{sec:particular}, we check that the main assumptions are not too strong in applying these results to particular set-valued maps $C(\cdot)$ given as the intersection of complements of convex sets. Moreover an application to a model of crowd motion is described. 

\gb This work was briefly presented in a proceeding \cite{Be}. 

\section{Preliminaries about Prox-regular sets} \label{sec:1}

We emphasize that the different notions defined in this section can be extended in the case of an infinite dimensional Hilbert space $H$. Here we only deal with a Euclidean framework. We denote by ${\mathbb B}$ the unit closed ball. For a subset $C$ of $\R^d$, we write $d_C$ for the distance function to this set:
$$d_C(x) := \inf_{y\in C} \ |y-x|.$$ 

\begin{df} Let $C$ be a closed subset of $\R^d$. The set-valued projection operator $P_C$ is defined on $\R^d$ by
$$ \forall x\in \R^d, \qquad P_C(x):=\left\{y\in C,\ |x-y|=d_C(x)\right\}.$$
\end{df}

\begin{df} Let $C$ be a closed subset of $\R^d$ and $x\in C$, we write $\NN(C,x)$ for the proximal normal cone of $C$ at $x$, defined by:
 $$ \NN(C,x):=\left\{v\in\R^d,\ \exists s>0, \ x\in P_C(x+sv)\right\}.$$
\end{df}

\mb We now come to the main notion of uniformly prox-regular set. It was initially introduced by H. Federer (in \cite{Federer}) in finite dimensional spaces under the name of ``positively reached set''. Then it was extended in infinite dimensional space and studied by F.H.~Clarke, R.J.~Stern and P.R.~Wolenski in \cite{Clarke} and by R.A. Poliquin, R.T.~Rockafellar and L.~Thibault in \cite{PRT}.

\begin{df} Let $C$ be a closed subset of $\R^d$ and $\eta>0$. The set $C$ is said $\eta$-prox-regular if for all $x\in C$ and $v\in \NN(C,x)\setminus \{0\}$
 $$B\left(x+\eta\frac{v}{|v|},\eta\right) \cap C = \emptyset.$$
Equivalently, $C$ is $\eta$-prox-regular if for all $y\in C$, $x\in \partial C$ and $v\in \NN(C,x)$
\be{eq:hypo} \langle y-x,v\rangle \leq \frac{|v|}{2\eta}|x-y|^2. \ee
\end{df}

\begin{rem} We refer the reader to \cite{Clarke,Clarke2} for other equivalent definitions related to the limiting normal cone. Moreover we can define the notion using the smoothness of the function distance $d_C$ (see \cite{PRT}).
\end{rem}
 This definition is very geometric, it describes the fact that we can continuously roll an external ball of radius smaller than $\eta$ on the whole boundary of the set $C$. The main property is the following one: for a $\eta$-prox-regular set $C$ and every $x$ satisfying $d_C(x)<\eta$, the projection of $x$ onto $C$ is well-defined (i.e. $P_C(x)$ is a singleton) and the projection is continuous.

\begin{df} \label{def} Let $A:\R^d \rightrightarrows \R^d$ be a set-valued operator on $\R^d$. We write $D(A)$ for its domain:
$$ D(A):=\left\{x\in\R^d,\ A(x) \neq\emptyset\right\}$$
and $\Gamma(A)$ for its graph:
$$ \Gamma(A):=\left\{(x,\alpha)\in\R^{2d},\ x\in D(A),\ \alpha\in A(x) \right\}.$$
The operator $A$ is said ``hypomonotone'' if there exists a constant $\delta>0$ such that for all $x,y\in D(A)$ and all $(\alpha,\beta)\in A(x)\times A(y)$, we have
\be{hypo}
\langle \alpha-\beta,x-y\rangle \geq -\delta\left[|\alpha|+|\beta|\right] |x-y|^2.
\ee
Such an operator is called hypomonotone with constant $\delta$ and it is maximal if for all $x,\alpha\in \R^d$ then
$$ \left( \forall y\in D(A),\ \forall \beta\in A(y),\  \langle \alpha-\beta,x-y\rangle \geq -\delta\left[|\alpha|+|\beta|\right] |x-y|^2 \right) \Longrightarrow \alpha\in A(x).$$
\end{df}

\mb Then  we refer the reader to the work \cite{PRT} of R.A. Poliquin, R.T. Rockafellar and L. Thibault for the following result:

\begin{prop} \label{prop:maxhypo} Let $C$ be a closed subset. Then $C$ is $\eta$-prox-regular if and only if the proximal normal cone $\NN(C,\cdot)$ is a maximal hypomonotone operator with constant $\frac{1}{2\eta}$. Moreover $D(\NN(C,\cdot))=C$.
\end{prop}

\mb The following definition comes from the work of A.A.~ Vladimirov (Section 3 in \cite{Vladimirov}). We extend it to hypomonotone set-valued operators.

\begin{df} \label{dist}  For $A$ and $B$ two hypomonotone operators with constant $\delta$, we define
$$d_{V}(A,B) := \sup_{(x,\alpha)\in \Gamma(A)} \ \sup_{(y,\beta)\in \Gamma(B)} \frac{\langle \alpha-\beta,y-x\rangle - 2\delta\left(|\alpha|+|\beta|\right)|x-y|^2}{1+|\alpha|+|\beta|} .$$
\end{df}

\mb Note that $d_V$ is not a distance as the triangle inequality is not satisfied. However, we have the following inequality:

\begin{prop} \label{p:ex} Let $C_1$ and $C_2$ be two $\eta$-prox-regular sets of $\R^d$ and write $d_H$ the Hausdorff distance
$$ d_H(C_1,C_2) := \sup_{y\in \R^d} \left|d_{C_1}(y) - d_{C_2}(y) \right|.$$
Then with $\delta=\frac{1}{2\eta}$, we have
$$ d_{V}(\NN(C_1,\cdot),\NN(C_2,\cdot)) \leq d_H(C_1,C_2)+\frac{1}{\eta} d_H(C_1,C_2)^2.$$
\end{prop}

\dem The proof is the same one as for Lemma 3.4 in \cite{Vladimirov} (dealing with convex sets). We detail it for an easy reference. Let $(x,\alpha)\in \Gamma(\NN(C_1,\cdot))$ and 
$(y,\beta)\in \Gamma(\NN(C_2,\cdot))$. We know that there exists a point $\tilde{x}\in P_{C_1}(y)\in C_1$ such that $|y-\tilde{x}|\leq  d_H(C_1,C_2)$. Then the hypomonotonicity property of the proximal normal cone implies that
\begin{align*} 
\langle \alpha , \tilde{x}-x \rangle & \leq \frac{1}{2\eta}|\alpha| |x-\tilde{x}|^2 \\
 & \leq \frac{1}{2\eta} |\alpha| \left[|x-y|+|y-\tilde{x}|\right]^2 \\
 & \leq \frac{1}{\eta} |\alpha| \left[|x-y|^2+|y-\tilde{x}|^2\right],
\end{align*}
hence
$$ \langle \alpha , y-x \rangle \leq |\alpha|d_H(C_1,C_2) + \frac{1}{\eta} |\alpha| |x-y|^2+\frac{1}{\eta} |\alpha|d_H(C_1,C_2)^2. $$
Since an analoguous inequality is satisfied by $\beta$, the desired result is obtained by summing both inequalities.
\findem

\begin{cor} Let $C:\I \rightrightarrows \R^d$ be a set-valued map taking $\eta$-prox-regular values. If $t\to C(t)$ is absolutely continuous (for the Hausdorff distance), then the operator-valued map $t\to \NN(C(t),\cdot)$ is absolutely continuous for $d_V$. 
\end{cor}

\mb According to the work of A.A.~ Vladimirov (\cite{Vladimirov}), an absolute continuity of the operator-valued map $t\to\NN(C(t),\cdot)$ is the appropriate assumption to solve the differential inclusion:
$$ \frac{dx}{dt}(t) + \NN(C(t),x(t)) \ni f(t,x(t)).$$
That is why, the previous proposition suggests to suppose an absolute continuity for the map $t\to C(t)$ (which is exactly the same assumption as done in \cite{Thibrelax,Thibbv}).

\mb We finish this section with the following result:

\begin{prop} \label{prop:regu}
Let $C$ be a $\eta$-prox-regular set in $\R^d$. Then for all $\epsilon< \eta/8$, the set
$$C_\epsilon:=\left\{x\in\R^d,\ d_C(x)\leq \epsilon \right\}=C+\epsilon{\mathbb B}$$
is uniformly $\eta/8$-prox-regular. Moreover the distance function:
$ d_{C_\epsilon}$ is $C^1$ in a neighborhood of $\partial C_\epsilon$. So the boundary $\partial C_\epsilon$ is a $C^1$-manifold.
\end{prop}

\begin{rem} This proposition is based on general results due to F. Bernard, L. Thibault and S. Zlateva in \cite{BTZ} and so it can be extended in a uniformly convex Banach space framework.
\end{rem}

\dem  We will use the notations of \cite{BTZ}, mainly for a closed set $C$, we define
$$ E_C(l):=\left\{ x\in\R^d,\ d_C(x)\geq l\right\}.$$
>From Theorem 6.2 and Lemma 6.9 of \cite{BTZ}, we deduce that for $l= \eta/4$ the set $E_{C}(l)$ is $l$-prox-regular. 
Then Lemma \ref{lem:EE} (later proved) ensures that for all $\epsilon\in(0,l/2)$ we have
$$ C_\epsilon=C+\epsilon{\mathbb B} = E_{{E_C}(l)}(l-\epsilon).$$
Analogously, we deduce that the set $C_\epsilon$ is $(l-\epsilon)$ prox-regular and so is $\eta/8$-prox-regular. \\
Moreover it can be checked that these sets have a $C^1$ boundary. Indeed Lemma \ref{lem:EE} yields
$$\partial C_\epsilon = \left\{x\in\R^d,\quad d_C(x)-\epsilon = 0 \right\}.$$
Applying Theorem 6.2 of \cite{BTZ}, we deduce that $d_C$ is $C^1$ in a neighbourhood of $\partial C_\epsilon$ and so $\partial C_\epsilon$ is a $C^1$ manifold. \findem

\mb The proof is ended provided that we prove the following lemma:

\begin{lem} \label{lem:EE} Let $C$ be a $\eta$-prox-regular set and $l=\eta/4$. Then for all $\epsilon\in(0,l/2)$
\be{eq:lem1} C_\epsilon:= \left\{x\in \R^d, \ d_C(x)\leq \epsilon \right\} = E_{E_{C}(l)}(l-\epsilon). \ee
Moreover 
\be{eq:lem2} \partial C_\epsilon:=\left\{x\in\R^d,\ d_C(x)= \epsilon \right\}. \ee
\end{lem}

\dem We first check the two embeddings of (\ref{eq:lem1}). First as the function $d_C$ is $1$-Lipschitz then for $\epsilon\leq l$, it yields
$$ \forall x\in C+\epsilon{\mathbb B},\ \forall u\in  E_{C}(l), \qquad l-\epsilon \leq d_C(u)-d_C(x) \leq |u-x|.$$
So we deduce that $d_{E_{C}(l)}(x)\geq l-\epsilon$ for all $x\in C+\epsilon{\mathbb B}$, which gives
$$C_\epsilon=C+\epsilon{\mathbb B} \subset E_{E_{C}(l)}(l-\epsilon).$$
For the reverse inclusion, we have to use the prox-regularity assumption. Let us take $x\in E_{E_{C}(l)}(l-\epsilon)$. We may suppose that $x$ does not belong to $C$, else $x\in C \subset C+\epsilon{\mathbb B}$. So let us consider $x_0\in \PPP_C(x)$. Then $x-x_0$ is a proximal normal vector at $x_0$ and so we know that
$$ u:= x_0 + l\frac{x-x_0}{|x-x_0|} \in E_{C}(l).$$
Here we have used that $l<\eta/2$ to get $d_C(u)=l$. Then as $x\in E_{E_{C}(l)}(l-\epsilon)$, we deduce that
$$ l-\epsilon \leq |u-x| = \left|l-|x-x_0|\right| = l-|x-x_0|,$$ 
where the last equality comes from the fact that $x$ does not belong to $E_C(l)$.
The inequality yields $|x-x_0|\leq \epsilon$, which ends the proof of
$$C+\epsilon{\mathbb B} \supset E_{E_{C}(l)}(l-\epsilon).$$
Finally, it remains us to prove (\ref{eq:lem2}). First thanks to the continuity of the distance function, it is obvious that
$$ \partial C_\epsilon \subset \left\{x\in\R^d,\ d_C(x)= \epsilon \right\}.$$
The previous reasoning (based on the prox-regularity of the set $C$) yields the reverse embedding and also concludes the proof. \findem

\section{Well-posedness and Stability results for stochastic sweeping process} \label{sec:2}

In all this section, we consider a set-valued map $C$. We first define the required assumptions, under which  well-posedness results for stochastic sweeping process can be proved. We want to study the stochastic differential inclusion (\ref{eq:sto}) which can be seen as the equation of a reflected Brownian motion onto the moving set $C(\cdot)$. In order to apply the results of P. L. Lions and A. S. Sznitman \cite{LS}, the sets $C(t)$ are supposed to satisfy some properties.

\mb We refer the reader to \cite{LS} for the following definitions without considering the time variable. Here we add the time-dependence.

\begin{df} \label{def:adm} The set-valued map $C$ is said {\em admissible} on $[0,T]$ if it takes uniformly prox-regular values (with a same constant) and if there exist $\delta,r,\tau>0$, and for all $t\in[0,T]$ sequences $(x_p)_p$ and $(u_p)_p$ with $|u_p|=1$ and $x_p\in C(t)$ such that for all $s\in [0,T]$ with $|t-s|\leq \tau$, $(B(x_p,r))_p$ is a bounded covering of the boundary $\partial C(s)$ and 
\be{hyp5} \forall p,\ \forall y\in \partial C(s) \cap B(x_p,2r),\ \forall v\in \NN(C(s),y), \quad \langle v,u_p\rangle \geq \delta |v|.\ee
\end{df}

\begin{rem} The original definition in \cite{LS} makes appear an extra assumption: there exists a sequence of approaching ``smooth sets'' satisfying a uniform bound of prox-regularity. As explained in \cite{Saisho}, the existence of such approaching smooth sets is not really necessary to prove the solvability of the Skorohod Problem. Moreover, due to the recent works about prox-regular sets, we know that such an approaching sequence always exists in a very general framework (see Proposition \ref{prop:regu}).
\end{rem}

\mb The second important property (used in \cite{LS}) is the following one:

\begin{df} \label{def:reg} A set-valued map $C$ is said {\em regular} on $I$ if there exists a function $\Phi\in C^2_b(I\times \R^d)$ satisfying for all $t\in[0,T]$, $x\in \partial C(t)$ and $v\in \NN(C(t),x)$
$$ \langle \nabla_x\, \Phi(t,x),v \rangle \leq -\mu |v|,$$ with $\mu>0$.
\end{df}

\mb According to Remark 3.1 of \cite{LS} and Lemma 5.3 of \cite{Saisho}, we know that an admissible set is locally regular. 
Moreover, Y. Saisho has proved in \cite{Saisho} that the local regularity is sufficient to obtain well-posedness results.

\begin{df}[\cite{Thibrelax, Thibbv}] \label{def:abs} A set-valued map $C(\cdot)$ is said to vary in an absolutely continuous way, if there exists an absolutely continuous function $v$ such that for all $t,s\in \I$
$$ d_H(C(t),C(s))\leq \left|v(t)-v(s)\right|,$$
where $d_H$ is the Hausdorff distance.
\end{df}

\subsection{Deterministic Skorohod problem: an extension of sweeping process} \label{subsec:21}

Before solving the stochastic differential inclusion (\ref{eq:sto}), we study the associated Skorohod problem. Let $h:\I\to \R^d$ be a continuous function satisfying $h(0)\in C(0)$. We say that a couple of continuous functions $(x,k)$ on $\I$ is a solution of the Skorohod problem (SkP,$h$) if:
 \begin{itemize} 
\item[$\bullet$] for all $t\in \I$, $x(t)\in C(t)$ 
\item[$\bullet$] the function $k$ is continuous and have a bounded variation on $\I$
\item[$\bullet$] the differential measure $dk$ is supported on $\left\{t\in\I,\ x(t)\in\partial C(t) \right\}$:
\be{eq1}  |k|(t) = \int_0^t {\bf 1}_{x(s)\in \partial C(s)} d|k|(s), \qquad k(t)=\int_0^t \xi(s) d|k|(s),\ee
with $\xi(s)\in \NN(C(s),x(s))$
\item[$\bullet$] for all $t\in \I$, we have
\be{eq2} x(t)+k(t)=h(t). \ee
\end{itemize} 

\mb Here we denote by $|k|(t)$ the total variation of the function $k$ on $[0,t]$. By extension, a continuous function $x$ is said to be a solution of (SkP,$h$) if there exists a function $k$ such that $(x,k)$ satisfies the previous properties.

\mb This subsection is devoted to the study of the deterministic problem (SkP,$h$) (defined by (\ref{eq1}) and (\ref{eq2})). By following the ideas of \cite{LS,Saisho}, we also begin with the following proposition:

\begin{prop} \label{prop3} Consider an admissible set-valued map $C$ varying in an absolutely continuous way and assume that for all $h\in C^\infty(\I,\R^d)$ with $h(0)\in C(0)$, there exists a solution $(x,k)$ to the Skorohod problem (SkP,$h$). \\
Then for all $h\in C^0(\I,\R^d)$ with $h(0)\in C(0)$, there exists a unique solution $(x,k)$ to the Skorohod problem (SkP,$h$). Furthermore the mapping $(h \rightarrow x)$ from $C^0(\I,\R^d)$ into itself is H\"older continuous of order $\frac{1}{2}$ on compact sets.
\end{prop}

\dem We refer the reader to Theorem 1.1 of \cite{LS} for a detailed proof of such proposition in the case of a constant set $C$. It is based on the hypomonotonicity and the admissibility properties. \\
\noindent {\bf First step:} Uniqueness. \\
First the uniqueness is ``as usual'' a direct consequence of the hypomonotonicity property and Gronwall's Lemma. We have to be careful as we are working with a function $k$ which is only assumed to have a bounded variation. \\
We recall the following version of Gronwall's Lemma (due to R. Bellman \cite{Bellman}):
\begin{lem} \label{lem:gronwall}
 Let $\I:=[0,T]$ be a closed bounded interval. Let $u$ be a non-negative measurable function defined on $\I$ and let $\mu$ be a locally finite non-negative measure on $\I$. Assume that $u\in L^1(\I,d\mu)$ and that for all $t\in \I$
$$ u(t) \leq \int_{[0,t]} u(s) d\mu(s).$$
If the function $t \to \mu([0,t])$ is continuous on $\I$, then for all $t\in \I$
$$ u(t) =0.$$
\end{lem}

\noindent Let us take two solutions $x$ and $\tilde{x}$ (we associate the corresponding functions $k$ and $\tilde{k}$). We study the error term:
$$ z(t):=x(t)-\tilde{x}(t) \qquad e(t):=|z(t)|^2.$$
The difference function $z$ is also solution of the following differential equation (in the sense of time-measure):
$$dz(t)= -dk(t)+d\tilde{k}(t).$$
As $k$ and $\tilde{k}$ have a bounded variation, it comes
$$ de(t) = 2\langle z(t),dz(t)\rangle = 2 \langle z(t), -dk(t)+d\tilde{k}(t) \rangle.$$
>From the hypomonotonicity of the proximal normal cones $\NN(C(t),\cdot)$ (due to the $\eta$-prox-regularity property of the sets $C(t)$, see Proposition \ref{prop:maxhypo}), we deduce that
$$de(t) \leq \frac{e(t)}{2\eta}\left(|dk|(t)+|d\tilde{k}|(t)\right).$$
Then as $k$ and $\tilde{k}$ are assumed to have finite variation on $\I$ and $e(0)=0$, we deduce that $e=0$ thanks to Lemma \ref{lem:gronwall} (with $\mu=|dk|+|d\tilde{k}|$). The proof of uniqueness is also concluded. \\

\mb {\bf Second step:} Existence. \\
Let $h\in C^\infty(\I,\R^d)$ and $(x,k)$ a solution to the Skorohod Problem (SkP,$h$). We follow the reasoning, used in the proof of Lemmas 1.1 and 1.2 of \cite{LS}, in pointing out the difficulties raised by the time-dependence of $C$. \\
First, Lemma 1.1 (6) of \cite{LS} still holds as it only rests on the uniform prox-regularity of sets $C(t)$. We detail how Lemma 1.1 (7) of \cite{LS} should be modified to take into account the time-dependence of $C$.
Let $\tau,r$ be given by the admissibility property (Definition \ref{def:adm}). On $[0,\tau]$, we denote by $({\mathcal O}_{i})_{i}$ the sets ${\mathcal O}_{i}:=B(x_i,2r)\cap C(0)$ and ${\mathcal O}_0$ the following set
$$ {\mathcal O}_{0} = C(0) \setminus \left[ \bigcup_i B(x_i,r)\right].$$
We let $T_1:= \inf \{t\in[0,\tau], \ x(t) \notin {\mathcal O}_{i_0} \}$, where $i_0$ is such that $h(0)=x(0)\in {\mathcal O}_{i_0}$. Then either $x(T_1)\in{\mathcal O}_{0}$ and we set $i_1=0$, or $x(T_1)\in B(x_{i_1},r)$ for some $i_1$. In this way we construct by induction, $i_m$ and $T_{i_m}$ such that if $T_m<\tau$, $x(T_m) \in B(x_{i_m},r)$ or $x(T_m) \in {\mathcal O}_{0}$, $T_{m+1}:=\inf \{ t\in ]T_m,\tau], \ x(t) \notin {\mathcal O}_{i_m}\}$.
Using the admissibility properties, we obtain such estimations for the variation of $k$: let $s,t\in[T_{m},T_{m+1}]$, $s\leq t$ for some $m$)
$$ |k|(t)-|k|(s) \leq \frac{1}{\delta}\left(|h(t)-h(s)| + |x(t)-x(s)|\right).$$
Moreover if $t_m=0$ then for all $t\in [T_m,T_{m+1}]$ $x(t)\notin \partial C(t)$ and the variation $|k|$ is constant on $[T_m,T_{m+1}]$ else for all $t\in [T_m,T_{m+1}]$ $x(t)\in B(x_{i_m},2r)$ and $|x(t)-x(s)|\leq 2r$.
We finally conclude that
\be{eq:K} |k|(T_{m+1})-|k|(T_m) \leq K, \ee
where $K$ is a numerical constant depending on $|h(T_{m+1})-h(T_{m})|$.
Then we can repeat the proof of Lemma 1.2 in \cite{LS} in order to give an upper bound to the variation $|k|(\tau)$. 
However the time-dependence of $C$ makes appear a new quantity in (7) of \cite{LS}. More precisely, we get for $s,t \in [T_m,T_{m+1}]$ \footnote{For two quantities $A,B$, we write $A\lesssim B$ if there exists a constant $C$ such that $A\leq CB$.}
\begin{align}
 |x(t)-x(s)|^2 \lesssim & \, \sup_{[s,t]} |h(\cdot)-h(s)|^2 + \sup_{[s,t]} |h(\cdot)-h(s)| \nonumber \\
 & + (|k|(t)-|k|(s))\left[\sup_{[s,t]} d_H(C(\cdot),C(s)) + \sup_{[s,t]} d_H(C(\cdot),C(s))^2\right]. \label{eq:newterm} 
\end{align}
Indeed, in the proof, we have to estimate for $u\in[s,t]$
$$ 2\langle x(u)-x(s),dk(u)\rangle - \eta d|k|(u) |x(u)-x(s)|^2.$$
To use the hypomonotonicity of the proximal normal cone (see (\ref{eq:hypo})), we have to write $x(s)=y(u)+z(u)$ with $y(u)\in \PPP_{C(u)}(x(s))$. As $x(s)\in C(s)$, $|z(u)|\leq d_H(C(u),C(s))$ and this operation makes appear a rest which is bounded by $d_H(C(u),C(s)) d|k|(u)$. \\
The two last terms in (\ref{eq:newterm}) can be assumed as small as we want since $C$ varies in an absolutely continuous way
and the variation of $k$ is bounded by (\ref{eq:K}). If $\tau$ is small enough (with respect to $r$), then by definition of $T_{m+1}$ and $T_m$ we deduce that for some constant $c$
$$ r^2 \leq |x(T_{m+1})-x(T_m)|^2 \leq c \sup_{[T_{m},T_{m+1}]} |h(\cdot)-h(s)|^2 + c\sup_{[T_m,T_{m+1}]} |h(\cdot)-h(s)| + r^2/2$$
which gives 
$$ r^2 \lesssim \sup_{[T_{m},T_{m+1}]} |h(\cdot)-h(s)|^2 + \sup_{[T_m,T_{m+1}]} |h(\cdot)-h(s)|.  $$
As done in \cite{LS}, we can also conclude that the length $|T_{m+1}-T_m|$ is not too small, according to the uniform continuity modulus of $h$ and so that the collection of indices $m$ is finite and its cardinal is bounded by a constant depending on $h$. Finally, the total variation $|k|_{[0,\tau]}:=|k|(\tau)-|k|(0)=|k|(\tau)$ is controlled by $\|h\|_\infty$ and the uniform continuity modulus of $h$. \\ The same reasoning holds on $[\tau,2\tau]$. Let us remark that if $(x,k)$ is the solution of (SkP,$h$) on $[0,2\tau]$ then $(x,k-k(\tau))$ is the solution of (SkP,$h-k(\tau)$) on $[\tau,2\tau]$. So we know that the variation $|k-k(\tau)|_{[\tau,2\tau]}=|k|_{[\tau,2\tau]}$ can be bounded by the $L^\infty$ norm of $h-k(\tau)$ (which is controlled by the $L^\infty$ norm and the continuity modulus of $h$) and its uniform continuity modulus (which is the same as for $h$). Consequently, the total variation $|k|_{[0,2\tau]}$ is bounded with the help of $h\in C^0([0,2\tau])$. By iterating (at most $T/\tau$ steps), we conclude that for $h\in C^\infty$, the (assumed) solution $(x,k)$ satisfies 
$$ \|x\|_{L^\infty([0,T])} + |k|_{[0,T]} \leq \tilde{K}, $$
where $\tilde{K}$ only depends on the $L^\infty$ norm and the uniform continuity modulus of $h$. This key-point permits us to take the limit and to prove Proposition \ref{prop3} by density arguments as in \cite{LS} (the proof can be easily adapted with a time-dependent set $C$). \findem


\begin{thm} \label{thm1} Consider an admissible set-valued map $C(\cdot)$, varying in an absolutely continuous way. For all $h\in C^0(\I,\R^d)$ with $h(0)\in C(0)$, there exists a unique solution $(x,k)$ to the Skorohod problem (SkP,$h$). Furthermore the mapping $(h \rightarrow x)$ from $C^0(\I,\R^d)$ into itself is H\"older continuous of order $\frac{1}{2}$ on compact sets. In addition if $h\in BV(\I)$ then $x\in BV(\I)$ and
\be{eq:aze} d|k|_t \lesssim d|h|_t + d|v|_t, \ee
where $v$ is given by Definition \ref{def:abs}.
\end{thm}

\dem According to the previous proposition, we just have to deal with smooth functions $h$. \\
Indeed for $h\in W^{1,1}(\I)$, the Skorohod problem (SkP,$h$) has already been treated since it corresponds to the so-called sweeping process. We refer the reader to \cite{Thibrelax,Thibbv} for more details about these differential inclusions. We deal with Lipschitz moving sets $C(t)$. 
So applying Theorem 1 of \cite{Thibrelax}, we know that there exists a pair of function $(x,k)\in W^{1,1}(\I)^2$ such that almost everywhere on $\I$
$$ \frac{dx}{dt} + \frac{dk}{dt} = \frac{dh}{dt}$$
and for almost every $t\in\I$
$$ \frac{dk}{dt}(t) \in \NN(C(t),x(t)).$$
Moreover, (\ref{eq:aze}) holds. \\
This gives us a solution to the Skorohod Problem (SkP,$h$) for every smooth function $h\in W^{1,1}(\I)$. Then we conclude the proof of the theorem, thanks to Proposition \ref{prop3}.
\findem


\subsection{Sweeping process with a stochastic perturbation} \label{sec:3}

In this section, we consider the Euclidean space $\R^d$ (equipped with its Euclidean structure), a probability space $\espace$ endowed with a standard filtration $(\F_t)_{t>0}$ and a standard $\R$-valued Brownian motion $(B_t)_{t>0}$ associated to this filtration. We denote by ${\mathbb E}(X)$ the expectation of a random variable $X$, according to this probability space. \\
We fix a bounded time-interval $\I=[0,T]$ and denote by $\BB$ the Banach space of time-continuous $\F_t$-adapted process $X$ satisfying
$$ \|X\|_\BB := \E \left[ \sup_{t\in \I} |X_t|^4 \right] ^{1/4}<\infty.$$
Let $f:\I\times \R^d \to \R^d$ and $\sigma:\I\times \R^d\to \R^d$ be two maps \footnote{We have chosen an $\R^d$-valued function $\sigma$ with a real Brownian motion $(B_t)_{t>0}$. Indeed all the results and the proofs hold for a real function $\sigma$ with an $\R^d$-valued Brownian motion.}. We look for solving the following stochastic differential inclusion on $\I$:

\be{eqsto} \left\{ \begin{array}{l} 
                  dX_t + \NN(C(t),X_t) \ni f(t,X_t)dt + \sigma(t,X_t) dB_t \vsp \\
                  X_0=u_0 ,
                 \end{array} \right.\ee
where $u_0\in C(0)$ is a non-stochastic initial data. Let us first give a more precise sense to this differential stochastic inclusion. 

\begin{df} \label{def:solution}
A continuous process $(X_t)_{t\in \I}$ is a solution of (\ref{eqsto}) if there exists another process $(K_t)_{t\in \I}$ such that:
\begin{enumerate}
 \item[a)] $(X_t)_{t\in \I}$ is a $\R^d$-valued process taking values in $C(t)$ and $\F_t$-adapted with continuous sample paths;
 \item[b)] $(K_t)_{t\in \I}$ is a $\R^d$-valued process, $\F_t$-adapted whose sample paths are continuous and have a bounded variation on $\I$;
 \item[c)] the following stochastic differential equation is satisfied
\be{eq2sto} dX_t + dK_t = f(t,X_t)dt + \sigma(t,X_t) dB_t; \ee
 \item[d)] the initial condition is verified: $X_0=u_0$ $\PP$-a.e.;
 \item[e)] the process $dK_t$ is supported on $\left\{t,\ X_t\in\partial C(t)\right\}$:
\be{eqK2}  |K|_t = \int_0^t {\bf 1}_{X_s\in \partial C(s)} d|K|_s, \qquad K_t=\int_0^t \xi(s) d|K|_s,\ee
with $\xi(s)\in \NN(C(s),X_s)$.
\end{enumerate}
\end{df}

\mb The point $e)$ gives a precise sense to ``$dK_t \in \NN(C(t),X_t)$''.

\mb Using the stochastic integral, we can rewrite (\ref{eq2sto}) as follows: we are looking for processes $(X_t)_{t\in \I}$ and $(K_t)_{t\in \I}$ satisfying that there exists 
a measurable set $\Omega_0\subset \Omega$ of full measure such that for all $\omega\in \Omega_0$ and all $t\in \I$ we have:
\be{eq3sto}   X_t=u_0 + \int_0^t f(s,X_s)ds + \int_0^t \sigma(s,X_s) dB_s - K_t. \ee

\subsubsection{Well-posedness results for (\ref{eqsto})}

\mb Now we come to our main results. We first assume that the set-valued map $C$ is admissible and regular in order to follow the ideas of P.L. Lions and A.S. Snitman in \cite{LS} and obtain well-posedness results in the space $L^4(\Omega,L^\infty(\I))$. We then describe results without requiring the regularity assumption following the ideas of Y. Saisho in \cite{Saisho}.

\begin{thm} \label{thm:principal} Consider an admissible and regular set-valued map $C$, varying in an absolutely continuous way. \\
Let $f,\sigma:\I\times \R^d \rightarrow \R^d$ be bounded and Lipschitz with respect to the second variable: there exists a constant $L$ such that for all $t\in \I$ and $x,y\in\R^d$
$$ \left|\sigma(t,x)-\sigma(t,y) \right| + \left|f(t,x)-f(t,y) \right| \leq L |x-y|$$
and
$$ \left|f(t,x)\right| + \left|\sigma(t,x)\right| \leq L.$$
Then  (\ref{eqsto}) is well-posed in $L^4(\Omega,L^\infty(\I))$. That means: for all initial data $u_0\in C(0)$, there exists one and only one process $(X_t)_{t\in \I}\in L^4(\Omega,L^\infty(\I))$ solution of (\ref{eqsto}) (in the sense of pathwise uniqueness). 
 \end{thm}


\mb {\bf Proof of Theorem \ref{thm:principal}:}
\mb 
For all $X\in\BB$, we write $F(X)$ the unique solution of the Skorohod Problem (SkP, $h$) with  
$$ F(X)_t +K_t = h(t):= u_0 +\int_0^t \sigma(s,X_s)dB_s + \int_0^t f(s,X_s)ds.$$ 
Due to Theorem \ref{thm1} and the continuity of the stochastic integral, the unique solution $F(X)$ exists. This well-defined map $F$ satisfies some properties as pointed out by the next proposition.

\begin{prop} \label{prop:imp} The map $F$ takes values in $\BB$.
There exists a constant $k$ such that for all processes $X,X'\in \BB$,
\be{contractante} \left\| F(X)-F(X') \right\|_\BB^4 \leq k \int_0^T \E \left[\|X-X'\|^4_{L^\infty([0,t])}\right] dt. \ee
The constant $k$ only depends on $|\I|=T$ and on the above constants about $C(\cdot)$, $f$ and $\sigma$.
\end{prop}

Consequently, for all integer $p\geq 2$, we have (by writing $F^{(p)}=F\circ .. \circ F$ for the iterated function):
\begin{align*}
 \left\| F^{(p)}(X)-F^{(p)}(X') \right\|_\BB^4 & \leq k \int_0^T \E\left[\|F^{(p-1)}(X)-F^{(p-1)}(X')\|^4_{L^\infty([0,t_1])}\right] dt_1 \\
 & \leq k^2 \int_0^T \int_0^{t_1} \E\left[\|F^{(p-2)}(X)-F^{(p-2)}(X')\|^4_{L^\infty([0,t_2])}\right] dt_2 dt_1 \\
 & \leq k^{p} \int_0^T \int_0^{t_1} \cdots \int_0^{t_{p-1}} \E\left[\|X-X'\|^4_{L^\infty([0,t_p])}\right] dt_p \cdots dt_1 \\
 & \leq \frac{T^pk^{p}}{p!} \|X-X'\|_{\BB}^4. 
\end{align*}
We also deduce that for a large enough integer $p$ (only dependent on the above constants), the map $F^{(p)}$ is a contraction on $\BB$. By the well-known fixed-point Theorem, we also know that $F$ admits one and only one fixed point $X\in \BB$. It is obvious that $X$ is then a solution of (\ref{eq3sto}) and (\ref{eqsto}). So we conclude the proof of Theorem \ref{thm:principal} provided that we prove Proposition \ref{prop:imp}. \findem

\mb {\bf Proof of Proposition \ref{prop:imp}: } \\
We only deal with (\ref{contractante}), a very similar reasoning permits to check that $F$ takes values in $\BB$. \\
We follow the ideas and the estimates of \cite{LS}. Let $X,X'\in \BB$, we denote $Y:=F(X)$, $Y':=F(X')$, $K,K'$ the associated processes and $\Phi$ the function given by Definition \ref{def:reg}. Then It\={o}'s formula implies
\begin{align} 
\Phi(t,Y_t) & = \Phi(0,u_0) + \int_0^t \frac{\partial \phi}{\partial t}(s,Y_s) ds + \int_0^t \langle \nabla_x \Phi(s,Y_s), f(s,X_s)\rangle ds \nonumber \\
 & + \int_0^t \langle \sigma(s,X_s), \nabla_x \Phi(s,Y_s)\rangle dB_s - \int_0^t \langle \nabla_x \Phi(s,Y_s), dK_s \rangle \\
 & + \frac{1}{2} \int_0^t \langle \sigma(s,Y_s), H_x\phi(s,Y_s)\sigma(s,Y_s)\rangle ds. \label{ito1} 
\end{align}
We denote the Hessian matrix by $H_x$. So with $\dsp \alpha:= \frac{\mu}{\eta}$ it comes 
\begin{align*} 
\exp\left\{\frac{-1}{\alpha}(\Phi(t,Y_t)+\Phi(t,Y'_t))\right\} |Y_t-Y'_t|^2  =  I+II+III
\end{align*}
with
\begin{align}
I:= &  2\int_0^t \exp\left\{\frac{-1}{\alpha}(\Phi(s,Y_s)+\Phi(s,Y'_s)) \right\} \Big[ \langle Y_s-Y'_s, \sigma(s,X_s)-\sigma(s,X'_s)\rangle dB_s  \label{ito2} \\
 & \hspace{0cm} + \langle Y_s-Y'_s, f(s,X_s)ds -f(s,X'_s)ds - dK_s + dK'_s \rangle \Big] \nonumber \\
 & \hspace{0cm} + \int_0^t \exp\left\{\frac{-1}{\alpha}(\Phi(s,Y_s)+\Phi(s,Y'_s)) \right\} tr\left[(\sigma(s,X_s)-\sigma(s,X'_s))(\sigma(s,X_s)-\sigma(s,X'_s))^t \right]ds \nonumber
\end{align}
corresponding to the differentiation of the square quantity $|Y_t-Y'_t|^2$, 
\begin{align}
II:= & \hspace{0cm} -\frac{1}{\alpha} \int_0^t \exp\left\{\frac{-1}{\alpha}(\Phi(s,Y_s)+\Phi(s,Y'_s))\right\} |Y_s-Y'_s|^2 \Big[ 
\langle \nabla_x \Phi(s,Y_s), f(s,X_s)\rangle ds   \nonumber \\
& \hspace{0cm} + \langle \nabla_x \Phi(s,Y'_s), f(s,X'_s)\rangle ds + \langle \sigma(s,X_s), \nabla_x \Phi(s,Y_s)\rangle dB_s + \langle \sigma(s,X'_s), \nabla_x \Phi(s,Y'_s)\rangle dB_s \nonumber \\
& \hspace{0cm} \frac{1}{2}\langle \sigma(s,X_s), H_x\phi(s,Y_s)\sigma(s,X_s)\rangle ds + \frac{1}{2}\langle \sigma(s,X'_s), H_x\phi(s,Y'_s)\sigma(s,X'_s)\rangle ds \nonumber \\
& \hspace{0cm} -  \langle \nabla_x \Phi(s,Y_s), dK_s\rangle - \langle \nabla_x \Phi(s,Y'_s), dK'_s\rangle \nonumber \\
& \hspace{0cm} + \frac{\partial \Phi}{\partial t}(s,Y_s) ds + \frac{\partial \Phi}{\partial t}(s,Y'_s) ds \Big] \nonumber \\
& \hspace{0cm} + \frac{1}{2\alpha^2} \int_0^t \exp\left\{\frac{-1}{\alpha}(\Phi(s,Y_s)+\Phi(s,Y'_s))\right\}|Y_s-Y'_s|^2  \label{var2} \\
 & \hspace{0.5cm} \Big[(\nabla_x\Phi(s,Y_s)\sigma(s,X_s)+\nabla_x\Phi(s,Y'_s)\sigma(s,X'_s))^2  \Big]ds \nonumber 
\end{align}
corresponding to the differentiation of the exponential quantity (with the It\={o}'s additional term in (\ref{var2})) and 
\begin{align}
III:= & \hspace{0cm} \frac{-2}{\alpha} \int_0^t \exp\left\{\frac{-1}{\alpha}(\Phi(s,Y_s)+\Phi(s,Y'_s))\right\} \label{var3}\\ 
& \hspace{0cm} \langle Y_s-Y'_s, \sigma(s,X_s)-\sigma(s,X'_s)\rangle \left[\nabla_x\Phi(s,Y_s)\sigma(s,X_s)+\nabla_x\Phi(s,Y'_s)\sigma(s,X'_s)\right] ds \nonumber
\end{align}
corresponding to the variational quadratic term due to the product between the square quantity and the exponential one.
Then we recall that $\Phi\in C^2_b$, by Definition \ref{def:reg} and that $C(t)$ is $\eta$-prox-regular, we have
\begin{align*}
 \frac{1}{\alpha}  \langle \nabla_x \Phi(s,Y_s), dK_s\rangle |Y_s-Y'_s|^2 - \langle Y_s-Y'_s,dK_s \rangle \leq 0 \\
 \frac{1}{\alpha}  \langle \nabla_x \Phi(s,Y'_s), dK'_s\rangle |Y_s-Y'_s|^2 - \langle Y_s-Y'_s,dK'_s \rangle \leq 0
\end{align*}
in the sense of a nonnegative time-measure.
So from the boundedness and the Lipschitz regularity of $\sigma$ and $f$, we deduce that 
\begin{align*}
|Y_t-Y'_t|^2  \lesssim &  \left|\int_0^t \exp\left\{\frac{-1}{\alpha}(\Phi(s,Y_s)+\Phi(s,Y'_s)) \right\} \langle Y_s-Y'_s, \sigma(s,X_s)-\sigma(s,X'_s)\rangle dB_s\right|    \\
  &  + \left| \int_0^t W_s |Y_s-Y'_s|^2 \left[\langle \sigma(s,X_s), \nabla_x \Phi(s,Y_s)\rangle + \langle \sigma(s,X'_s), \nabla_x \Phi(s,Y'_s)\rangle\right]dB_s\right| \\
  &  + \int_0^t |Y_s-Y'_s| |X_s-X'_s| ds + \int_0^t | Y_s-Y'_s|^2 ds  \\
  &  + \int_0^t |X_s-X'_s|^2 ds. 
\end{align*}
We denote by $W_s$
$$W_s:=\exp\left\{\frac{-1}{\alpha}(\Phi(s,Y_s)+\Phi(s,Y'_s))\right\}$$. 
\begin{align*}
\E\left[\sup_{t\leq t_0} |Y_t-Y'_t|^4\right]  & \lesssim  \E \left[\sup_{t\leq t_0}\left|\int_0^t W_s \langle Y_s-Y'_s, \sigma(s,X_s)-\sigma(s,X'_s)\rangle dB_s\right|^2\right]    \\
 &  \ \ \ +  \E \left[\sup_{t\leq t_0}\left| \int_0^t W_s |Y_s-Y'_s|^2 \langle \sigma(s,X_s), \nabla_x \Phi(s,Y_s)\rangle dB_s \right|^2\right] \\
 &  \ \ \ +  \E \left[\sup_{t\leq t_0}\left| \int_0^t W_s |Y_s-Y'_s|^2 \langle \sigma(s,X'_s), \nabla_x \Phi(s,Y'_s) dB_s \right|^2\right] \\
 &  \ \ \  +  \int_0^{t_0} \E\left[\sup_{t\leq s} | Y_t-Y'_t|^4 \right] ds  + \int_0^{t_0} \E\left[\sup_{t\leq s} |Y_t-Y'_t|^2 |X_t-X'_t|^2 \right] ds \\
 &  \ \ \  +  \int_0^{t_0} \E\left[\sup_{t\leq s} | X_t-X'_t|^4 \right] ds 
\end{align*}
Then by noting that the three first terms are submartingales, we can apply Doob's inequality in order to obtain (with Cauchy-Schwartz inequality) for every $t_0\in \I$:
\begin{align*}
\E\left[\sup_{t\leq t_0} |Y_t-Y'_t|^4\right]  & \lesssim  \E \left[\sup_{t\leq t_0}\left|\int_0^t \langle Y_s-Y'_s, \sigma(s,X_s)-\sigma(s,X'_s)\rangle dB_s\right|^2\right]    \\
 &  \ \ \  +  \int_0^{t_0} \E\left[\sup_{t\leq s} | Y_t-Y'_t|^4 \right] ds  + \int_0^{t_0} \E\left[\sup_{t\leq s} |Y_t-Y'_t|^2 |X_t-X'_t|^2 \right] ds \\
 &  \ \ \  +  \int_0^{t_0} \E\left[\sup_{t\leq s} | X_t-X'_t|^4 \right] ds \\
 & \lesssim   \E \left[\int_0^{t_0} \left|Y_s-Y'_s\right|^2|X_s-X_s'|^2 ds\right]    \\
 &  \ \ \  +  \int_0^{t_0} \E\left[\sup_{t\leq s} | Y_t-Y'_t|^4 \right] ds  + \int_0^{t_0} \E\left[\sup_{t\leq s} |Y_t,Y'_t|^2 |X_t-X'_t|^2 \right] ds \\
 &  \ \ \  +  \int_0^{t_0} \E\left[\sup_{t\leq s} | X_t-X'_t|^4 \right] ds. 
\end{align*}
The inequality $ xy\leq  2y^2 + 2x^2$ implies
$$ \E\left[\sup_{t\leq t_0} |Y_t-Y'_t|^4\right] \lesssim  \int_0^{t_0}  \E \left[\sup_{t\leq s}\left|Y_t-Y'_t\right|^4 \right] ds + \int_0^{t_0}  \E \left[\sup_{t\leq s}\left|X_t-X'_t\right|^4 \right] ds $$ 
which with Gronwall's Lemma, gives us
$$ \E\left[\sup_{t\leq t_0} |Y_t-Y'_t|^4\right]  \lesssim \int_0^{t_0}  \E \left[\sup_{t\leq s}\left|X_t-X'_t\right|^4 \right] ds. $$ 
This concludes the proof of Proposition \ref{prop:imp}. \findem

\mb We finish this subsection by results without requiring the regularity of the set-valued map $C$.

\begin{thm} \label{thm:principal2} Consider an admissible set-valued map $C$, varying in an absolutely continuous way. \\
Let $f,\sigma:\I\times \R^d \rightarrow \R^d$ be bounded and Lipschitz with respect to the second variable: there exists a constant $L$ such that for all $t\in \I$ and $x,y\in\R^d$
$$ \left|\sigma(t,x)-\sigma(t,y) \right| + \left|f(t,x)-f(t,y) \right| \leq L |x-y|$$
and
$$ \left|f(t,x)\right| + \left|\sigma(t,x)\right| \leq L.$$
Then (\ref{eqsto}) is well-posed. That means: for all initial data $u_0\in C(0)$, there exists one and only one process $(X_t)_{t\in \I}$ solution of (\ref{eqsto}) (in the sense of pathwise uniqueness). 
 \end{thm}

\dem Theorem \ref{thm:principal2} is a consequence of Theorem \ref{thm1}. Indeed it suffices to repeat the proof developed in \cite{Saisho}, which permits to obtain Theorem 5.1 \cite{Saisho} (with abstract compactness arguments) as a consequence of Theorem 4.1 \cite{Saisho} (which corresponds to the well-posedness of deterministic Skorohod problem).
\findem

\begin{rem} As already pointed out in Remark 5.1 \cite{Saisho}, the existence part of Theorem \ref{thm:principal2} only requires the boundedness and the continuity of functions $f,\sigma$. The Lipschitz regularity is necessary only for the uniqueness. 
 \end{rem}

\subsubsection{Stability results}

\mb Now we are looking for stability results: we let the stochastic perturbation $\sigma$ goes to $0$ and we prove that the corresponding solution tends to the solution of the deterministic sweeping process.

\begin{thm} \label{thm:conv1} Let $\sigma^\epsilon:\I\times \R^d \rightarrow \R^d$ be maps satisfying there exists $L$ with for all $\epsilon>0$
$$ \left|\sigma^{\epsilon}(t,x)-\sigma^\epsilon(t,y) \right| \leq L |x-y|$$
and
$$ \left\|\sigma^\epsilon\right\|_{L^\infty(\I\times \R^d)} \leq L.$$
We assume that $\sigma^\epsilon$ tends to $0$ in $L^\infty(\I\times \R^d)$ when $\epsilon$ goes to $0$. 
Under the assumptions of Theorem \ref{thm:principal}, we consider a fixed initial data $u_0\in C(0)$. For each $\epsilon\in(0,1]$, we denote $(X^\epsilon_t)_{t\in \I}$ the (unique) process solution of
$$ \left\{ \begin{array}{l} 
                  dX^\epsilon_t + \NN(C(t),X^\epsilon_t) \ni f(t,X_t)dt + \sigma^\epsilon(t,X^\epsilon_t) dB_t \vsp \\
                  X^\epsilon_0=u_0 \ .
                 \end{array} \right. $$
Denote $x$ the solution of 
$$ \left\{ \begin{array}{l} 
                 \dsp \frac{dx}{dt}(t) + \NN(C(t),x(t)) \ni f(t,x(t)) \vsp \\
                  x(0)=u_0 \ ,
                 \end{array} \right. $$
given by Theorem 1.1 of \cite{Thibrelax} and consider the deterministic process: defined for all $\omega\in\Omega$ by $X_t(\omega)=x_t$. \\
Then $X^\epsilon$ converges to $X$ in $L^4(\Omega,L^\infty(\I))$:
\be{eq:thm:conv1} \left\| X^\epsilon - X\right\|_{L^4(\Omega,L^\infty(\I))} \leq c_{u_0} \|\sigma^\epsilon\|_{L^\infty(\I\times \R^d)} \xrightarrow[\epsilon\to 0]{} 0, \ee
for some constant $c_{u_0}$ independent on $\epsilon$.
\end{thm}

\dem We denote $F(\epsilon,\cdot)$ the map defined on ${\mathcal B}:=L^4(\Omega,L^\infty(\I))$ into ${\mathcal B}$ as follows:
$\left[F(\epsilon,X)\right]_t$ is the unique solution of the Skorohod Problem  
$$ F(\epsilon, X)_t +K^{\epsilon}_t = u_0 +\int_0^t \sigma^\epsilon(s,X_s)dB_s + \int_0^t f(s,X_s)ds.$$ 
Analogously we define $F(0,\cdot)$. \\
Then in the previous subsection, we have proved that $X^\epsilon$ is the unique fixed point of $F(\epsilon,\cdot)$ and similarly $X$ is the unique fixed point of $F(0,\cdot)$. 
In order to apply the ``fixed point theorem with parameter'', we check that the map $F$ is continuous at $\epsilon=0$. \\
By Proposition \ref{prop:imp}, there exists a large enough integer $p$ (only depending on the above constants and $L$) such that for all $\epsilon\geq 0$, the map $F(\epsilon,\cdot)^{(p)}$ is a $\frac{1}{2}$-Lipschitzean map on $\BB$. We also deduce that
\begin{align*}
 \left\|X^\epsilon-X\right\|_{\mathcal B} & = \left\| F^{(p)}(\epsilon,X^\epsilon) - F^{(p)}(0,X) \right\|_{\mathcal B} \\
 & \leq \left\| F^{(p)}(\epsilon,X^\epsilon) - F^{(p)}(\epsilon,X) \right\|_{\mathcal B} + \left\| F^{(p)}(\epsilon,X)-F^{(p)}(0,X) \right\|_{\mathcal B} \\
 &  \leq \frac{1}{2}\left\| X^\epsilon-X \right\|_{\mathcal B} + \left\| F^{(p)}(\epsilon,X)-F^{(p)}(0,X) \right\|_{\mathcal B},
\end{align*}
which gives
$$ \left\|X^\epsilon-X\right\|_{\mathcal B} \leq 2 \left\| F^{(p)}(\epsilon,X)-F^{(p)}(0,X) \right\|_{\mathcal B}. $$
It also remains to bound this quantity. By Proposition \ref{prop:imp} the sequence $(F(\epsilon,X))_\epsilon$ is uniformly bounded in $\BB$ and the map $F(\epsilon,\cdot)$ are uniformly Lipschitz in $\BB$ (with a constant denoted by $k$). Hence, Proposition \ref{prop:conv1} (see below) implies
\begin{align*}
\left\| F^{(2)}(\epsilon,X) - F^{(2)}(0,X) \right\|_{\mathcal B} & \\
 & \hspace{-3cm} \leq  \left\| F(\epsilon,F(\epsilon,X)) - F(\epsilon,F(0,X)) \right\|_{\mathcal B}+\left\| F(\epsilon,F(0,X)) - F(0,F(0,X)) \right\|_{\mathcal B} \\
 & \hspace{-3cm} \leq k\tilde{k}\|X\|_{\mathcal B} \|\sigma^\epsilon\|_{L^\infty(\I\times \R^d)} + \tilde{k}\|\sigma^\epsilon\|_{L^\infty(\I\times \R^d)} \|F(0,X)\|_{\mathcal B}.
\end{align*}
Then by iterating the reasoning, we deduce that there exists a constant $c_{u_0}$ (depending on $u_0$ through $X$, $F(0,X)$, ..., $F^{(p)}(0,X)$ and on $p$) such that (\ref{eq:thm:conv1}) holds.
\findem

\begin{prop} \label{prop:conv1} With the notations of Theorem \ref{thm:conv1}, there exists a constant $\tilde{k}$ (independent on $\epsilon$) such that
\be{eq:conv1} \left\| F(\epsilon,X)-F(0,X) \right\|_\BB\leq \tilde{k} \|X\|_\BB \|\sigma^\epsilon\|_{L^\infty(\I\times \R^d)}.\ee
\end{prop}

\dem Let us denote the processes $Y^\epsilon:=F(\epsilon,X)$, $Y:=F(0,X)$ and $Z:=F(\epsilon,X)-F(0,X)$. Then $Z$ satisfies the following stochastic differential equation:
$$ dZ_t = -\left(dK^\epsilon_t - dK^0_t\right) + \sigma^\epsilon(t,X_t)dB_t.$$
With the same arguments as in the proof of Proposition \ref{prop:imp}, it can be shown 
\begin{align*}
|Z_t|^2=|Y^\epsilon_t-Y_t|^2  \lesssim & \  \left|\int_0^t \langle Y^\epsilon_s-Y_s, \sigma^\epsilon(s,X_s)\rangle dB_s\right|    \\
  &  + \int_0^t | Y^\epsilon_s-Y_s|^2 ds + \|\sigma^\epsilon\|_{L^\infty(\I\times \R^d)}^2. 
\end{align*}
Then applying Doob's inequality, we get for every $t_0\in \I$
\begin{align*}
\E\left[\sup_{t\leq t_0} |Z_t|^4\right]  & \lesssim  \E \left[\sup_{t\leq t_0}\left|\int_0^t \langle Y^\epsilon_s-Y_s, \sigma^\epsilon(s,X_s)\rangle dB_s\right|^2\right]    \\
 &  \ \ \  +  \int_0^{t_0} \E\left[\sup_{t\leq s} | Z_t|^4 \right] ds +  \|\sigma^\epsilon\|_{L^\infty(\I\times \R^d)}^4  \\
 & \lesssim   \|\sigma^\epsilon\|_{L^\infty(\I\times \R^d)}^2 \E \left[\int_0^{t_0} \left|Z_s\right|^2 ds\right]    \\
 &  \ \ \  +  \int_0^{t_0} \E\left[\sup_{t\leq s} | Z_t|^4 \right] ds + \|\sigma^\epsilon\|_{L^\infty(\I\times \R^d)}^4.
\end{align*}
By Cauchy-Schwartz inequality, we obtain with another constant $c$ (depending on $|\I|$)
\begin{align*}
\E\left[\sup_{t\leq t_0} |Z_t|^4\right] \lesssim  \|\sigma^\epsilon\|_{L^\infty(\I\times \R^d)}^4 + \int_0^{t_0} \E\left[\sup_{t\leq s} | Z_t|^4 \right] ds.
\end{align*}
By Gronwall's lemma, we obtain (\ref{eq:conv1}). \findem

\section{Euler scheme for stochastic sweeping process} \label{sec:schema}

We refer the reader to the work of Y.~Saisho \cite{Saisho}, where the convergence of some discretized Skorohod problems to the continuous Skorohod problem is studied with a constant set $C$. For $l\in C^0(\I)$ and a partition of $\I$ given by $(t^{n}_h=nh)_{1\leq n\leq Th^{-1}}$, we denote by $x_h$ the following function:
$$ \left\{ \begin{array}{l}
            x_h(t) := u_0 ,\quad \textrm{ for all } t\in[0,t^{1}_h] \vsp \\
            x_h(t) : = \PPP_{C}\left[x_h(t^n_h) + l(t) - l(t^n_h) \right] ,\quad \textrm{ for all } t\in[t^{n}_h, t^{n+1}_h].
          \end{array}
\right.  
$$
Then in \cite{Saisho}, it is proved that $x_h$ strongly converges in $L^\infty(\I)$ to the unique solution $x$ of the Skorohod problem
\be{Sk} \left\{ \begin{array}{l}
            x(t) + k(t) = l(t), \quad t\in\I \vsp \\
            x(0)=u_0.
           \end{array} \right. \ee
Moreover, the associated map $k_h:=l-x_h$ strongly converges in $L^\infty(\I)$ to $k$. 

\mb This section is devoted to the extension of such results for time-dependent sets $C(\cdot)$. Let us consider a set-valued map $C(\cdot)$ on $\I=[0,T]$. We now define a discretized solution $\tilde{x}_h$ as follows:
\be{schema} \left\{ \begin{array}{l}
            \tilde{x}_h(t) := u_0 ,\quad \textrm{ for all } t\in[0,t^{1}_h] \vsp \\
            \tilde{x}_h(t) : = \PPP_{C(t^{n+1}_h)}\left[\tilde{x}_h(t^n_h) + l(t) - l(t^n_h) \right] ,\quad \textrm{ for all } t\in]t^{n}_h, t^{n+1}_h].
          \end{array}
\right.  
\ee
If the set-valued map $C$ varies in an absolutely continuous way and takes uniformly prox-regular values, then it can be checked that $\tilde{x}_h(t^n_h) + l(t) - l(t^n_h)$ is close to $C(t^{n+1}_h)$ and so its projection is single-valued for $h$ small enough. The scheme is also well-defined for $h$ small enough and semi-implicit as we consider $C(t^{n+1}_h)$. It is a prediction-correction algorithm: predicted point $\tilde{x}_h(t^n_h) + l(t) - l(t^n_h)$, that may not belong to $C(t^{n+1}_h)$, is projected onto $C(t^{n+1}_h)$. Moreover, we refer the reader to the works (mentioned in the introduction) dealing with deterministic sweeping process. Such schemes are well-known in the framework of sweeping process (when $l\in W^{1,1}(\I)$) and corresponds to the so-called {\em Catching-up Algorithm} introduced by J.J Moreau in \cite{Moreausweep}.

\begin{rem} \label{rem:rem2}
 If $C(\cdot)$ is a set-valued map taking uniformly prox-regular values and varying in an absolutely continuous way, it is well-known that the Euler scheme is convergent for smooth functions $l$ (as in this case, Skorohod problem corresponds to sweeping process), see for example \cite{Thibrelax,Thibbv}: for $l\in W^{1,1}(\I)$, the sequence $(\tilde{x}_h)_h$ strongly converges to the unique solution $x$ of (\ref{Sk}). 
\end{rem}

\begin{thm} \label{thmm} Consider an admissible set-valued map $C$, varying in an absolutely continuous way. For every function $l\in C^0(\I)$, the discretized solutions $(\tilde{x}_h)_h$ strongly converges in $L^\infty(\I)$ to the unique solution $x$ of (\ref{Sk}).
\end{thm}

\dem According to Proposition \ref{prop:123}, $(\tilde{x})_h$ is uniformly continuous with respect to $h$ on $\I$. By Arzela-Ascoli's Theorem, $\{\tilde{x}_h, h>0\}$ is relatively compact in $C^0(\I)$. Then by uniqueness of solution to the continuous problem (\ref{Sk}), it suffices to show that every convergent subsequence (still denoted by $(\tilde{x}_h)_{h}$) converges to a solution of (\ref{Sk}). Let us fix such a convergent subsequence and write $x$ for its limit. Obviously $x$ is continuous. Since for $t\in]t_h^n,t_h^{n+1}]$,
$ \tilde{x}_h(t) \in C(t_h^{n+1})$, we deduce that
$$ d_{C(t)}(\tilde{x}_h(t)) \leq d_H(C(t_h^{n+1}), C(t)) \leq \left|v(t_h^{n+1})- v(t)\right| \xrightarrow[h\to 0]{} 0,$$
where $v$ is the continuous function given by Definition \ref{def:abs}. So we get that for all $t\in \I$, $x(t)\in C(t)$. Moreover Proposition \ref{prop:123} implies that $k=l-x$ has a bounded variation on $\I$ so it suffices to check (\ref{eq1}). Indeed the algorithm implies that with $k_h:=l-\tilde{x}_h$ and $t\in]t^n_h,t^{n+1}_h]$
$$ k_h(t) - k_h(t^n_h) \in \NN(C(t_h^{n+1}),\tilde{x}_h(t)).$$
We let the details to the reader and we refer to Theorem 4.1 of \cite{Saisho} for arguments in order to take the limit in the previous discrete inclusion and to deduce (\ref{eq1}). The idea rests on the hypomonotonicity of the proximal normal cone. \\
This permits to prove that $x$ is solution of (\ref{Sk}) and so by compactness we deduce that $\tilde{x}_h$ converges to the unique solution of (\ref{Sk}).
 \findem

\begin{rem} We would like to describe another way than the compactness argument, allowing us to prove the previous Theorem. Indeed by a similar reasoning than the one used in Proposition \ref{prop:123}, we can obtain a discrete version of (6) Lemma 1.1 in \cite{LS}: let $l,\overline{l}$ two continuous functions, we denote $\overline{x}$ and $\tilde{\overline{x}}_h$, the continuous and discretized solutions of the Skorohod problem (SkP,$\overline{l}$), then
$$ \|\tilde{\overline{x}}_h-\tilde{x}_h\|_{L^\infty(\I)}^2 \lesssim \|l-\overline{l}\|_{L^\infty(\I)} + \|l-\overline{l}\|_{L^\infty(\I)}^2.$$
Moreover the implicit constant does not depend on $h$ and depends only on $l$ and $\overline{l}$. So we conclude that the maps $(l \rightarrow \tilde{x}_h)$ is H\"older continuous of order $\frac{1}{2}$ from $C^0(\I,\R^d)$ into itself on compact sets (and the bound can be chosen independent on $h$). Thanks to Proposition \ref{prop3}, the map $(l \rightarrow x)$ is H\"older continuous of order $\frac{1}{2}$ from $C^0(\I,\R^d)$ into itself on compact sets too. \\
So for $l\in C^0(\I)$ and $\epsilon>0$, there exists a smooth function $\overline{l}$ such that $\|l-\overline{l}\|_{L^\infty(\I)}\leq \epsilon$. Then we have
$$ \| \tilde{x}_h-x \|_{L^\infty(\I)} \leq \|\tilde{x}_h-\tilde{\overline{x}}_h\|_{L^\infty(\I)} + \| \tilde{\overline{x}}_h - \overline{x}\|_{L^\infty(\I)} + \|\overline{x} - x\|_{L^\infty(\I)}.$$
>From the H\"older regularity of the maps $(l \rightarrow \tilde{x}_h)$ and $(l \rightarrow x)$, we know that the first and third terms are bounded by $\epsilon^{1/2}$ (up to a numerical constant). Thanks to Remark \ref{rem:rem2}, the second term tends to $0$ as $\overline{l}\in W^{1,1}(\I)$. Thus we deduce that $\tilde{x}_h$ uniformly converges to $x$.
\end{rem}

\begin{prop} \label{prop:123} The continuous functions $\tilde{x}_h$ are uniformly continuous on $\I$ and $\tilde{x}_h-l$ has a bounded variation, both them uniformly with respect to $h$. These estimations only depend on $l$ via its $L^\infty$-norm  and its uniform continuity modulus. 
\end{prop}

\noindent This proposition can be seen as a ``discrete version'' of (\ref{eq:newterm}) (or Lemma 1.1 (7) and Lemma 1.2 of \cite{LS}), so we will follow their proofs.

\dem  We take again the notations, introduced in the proof of Proposition \ref{prop3} (second step). We set $k_h:=l-\tilde{x}_h$ and let $m$ be a fixed integer. \\
{\bf First step:} Boundedness of a discrete variation $|k_h|(T_{m+1})-|k_h|(T_m)$.\\
For every $s,t\in [T_m,T_{m+1}]$, thanks to the admissibility property, we get (with $n_0,n_1$ satisfying $t_h^{n_0}\in[s,s+h[$ and $t_h^{n_1}\in]t-h,t]$):
\begin{align*} 
|\tilde{x}_h(t)-\tilde{x}_h(s)| + |l(t)-l(s)| & \geq \langle \tilde{x}_h(s)-\tilde{x}_h(t)+l(t)-l(s), u_{i_m}\rangle \\
 & \geq \sum_{n=n_0}^{n_1-1} \langle -\tilde{x}_{h}(t^{n+1}_h)+\tilde{x}_{h}(t^n_h) + l(t^{n+1}_h)-l(t^n_h) , u_{i_m} \rangle \\
& \ + \langle -\tilde{x}_h(t)+\tilde{x}_{h}(t^{n_1}_h), u_{i_m}\rangle  + \langle -\tilde{x}_{h}(t_h^{n_0})+\tilde{x}_{h}(s),u_{i_m}\rangle \\
& \ + \langle l(t)-l(t^{n_1}_h), u_{i_m}\rangle + \langle l(t_h^{n_0})-l(s),u_{i_m}\rangle \\
& \geq  \delta \sum_{n=n_0}^{n_1-1} |\tilde{x}_{h}(t^{n+1}_h)-\tilde{x}_{h}(t^n_h) - l(t^{n+1}_h)+l(t^n_h)| \\
& \ +\delta|\tilde{x}_h(t)-\tilde{x}_{h}(t^{n_1}_h)-l(t)+l(t^{n_1}_h)|\\
& \ +\delta|\tilde{x}_{h}(t_h^{n_0})-\tilde{x}_{h}(s)-l(t_h^{n_0})+l(s)|.
\end{align*}
We have used that for all integer $n$ and $t\in]t_h^n,t_h^{n+1}]$
$$ \tilde{x}_{h}(t)-\tilde{x}_{h}(t_h^n)-l(t)+l(t_h^n) \in -\NN(C(t_h^{n+1}),\tilde{x}_h(t)).$$
So we conclude 
$$  |k_h|(t)-|k_h|(s) \leq \frac{1}{\delta} \left( |\tilde{x}_h(t)-\tilde{x}_h(s)| + |l(t)-l(s)| \right), $$
where we set the discrete variation
$$ |k_h|(t)-|k_h|(s) := |k_h(t)-k_h(t^{n_1}_h)| + \sum_{n=n_0}^{n_1-1} |k_h(t^{n+1}_h)-k_h(t^n_h)| + |k_h(s)-k_h(t^{n_0}_h)|.$$
Consequently,
\be{eq:var} |k_h|(T_{m+1}) - |k_h|(T_{m}) \leq K, \ee
for some numerical constant, as for Proposition \ref{prop3}. \\
{\bf Second step:} Uniform continuity of $\tilde{x}_h$.\\
Let $\eta$ be a constant of prox-regularity of all the sets $C(\cdot)$. Furthermore for $s,t\in [T_m,T_{m+1}]$
we write
$$ |\tilde{x}_h(t)-\tilde{x}_h(s)|^2e^{-[|k_h|(t)-|k_h|(s)]/\eta} =  |\tilde{x}_h(t^{n_1}_h)-\tilde{x}_h(t_h^{n_0})|^2 e^{-[|k_h|(t^{n_1}_h)-|k_h|(t^{n_0}_h)]/\eta} + \textrm{Rest}.$$
We denote $S$ for the first term and we only study it as the {\em Rest} can be similarly estimated. By a discrete differentiation, it comes
\begin{align*} 
S = & \sum_{n=n_0}^{n_1-1} \left(|\tilde{x}_h(t^{n+1}_h)-\tilde{x}_h(t^{n_0}_h)|^2 - |\tilde{x}_h(t^{n}_h)-\tilde{x}_h(t^{n_0}_h)|^2\right)e^{-[|k_h|(t^{n+1}_h)-|k_h|(t^{n_0}_h)]/\eta} \\
 & + \sum_{n=n_0}^{n_1-1}|\tilde{x}_h(t^{n}_h)-\tilde{x}_h(t^{n_0}_h)|^2 \left(e^{- [|k_h|(t^{n+1}_h)-|k_h|(t^{n_0}_h)]/\eta}-e^{-[|k_h|(t^{n}_h)-|k_h|(t^{n_0}_h)]/\eta}\right).
\end{align*}
Furthermore
\begin{align} 
S & =  \sum_{n=n_0}^{n_1-1} \left\langle \tilde{x}_h(t^{n+1}_h)-\tilde{x}_h(t^{n}_h), \tilde{x}_h(t^{n+1}_h)+\tilde{x}_h(t^{n}_h)-2\tilde{x}_h(t^{n_0}_h)\right\rangle e^{-[|k_h|(t^{n+1}_h)-|k_h|(t^{n_0}_h)]/\eta} \label{no1} \\
& \ + \sum_{n=n_0}^{n_1-1} - |\tilde{x}_h(t^{n}_h)-\tilde{x}_h(t^{n_0}_h)|^2 \left[e^{- [|k_h|(t^{n}_h)-|k_h|(t^{n_0}_h)]/\eta}-e^{-[|k_h|(t^{n+1}_h)-|k_h|(t^{n_0}_h)]/\eta}\right] \label{no2} \\
& := S_1+S_2. 
\end{align}
Then in $S_1$, we make the following replacement in the first term in the inner product
\be{no11} \tilde{x}_h(t^{n+1}_h)-\tilde{x}_h(t^{n}_h) = \left(l(t^{n+1}_h)-l(t^{n}_h)\right) - \left(k_h(t^{n+1}_h)-k_h(t^{n}_h)\right). \ee
Thus we write $S_1=I+II$.
The first term $I$ equals to
\begin{align*}
 I = & \sum_{n=n_0}^{n_1-1} \langle l(t^{n+1}_h) - l(t^n_h),l(t^{n+1}_h)+l(t^{n}_h)-2l(t^{n_0}_h) \rangle e^{- [|k_h|(t^{n+1}_h)-|k_h|(t^{n_0}_h)]/\eta} \\
& - \langle l(t^{n+1}_h) - l(t^n_h), k_h(t^{n+1}_h)+k_h(t^{n}_h)-2k_h(t^{n_0}_h) \rangle e^{- [|k_h|(t^{n+1}_h)-|k_h|(t^{n_0}_h)]/\eta}.
\end{align*}
The first quantity $I_1$ (by producing the reverse manipulations with $l$ instead of $\tilde{x}_h$) can be bounded by $\sup_{s\leq t_1\leq t_2 \leq t} |l(t_2)-l(t_1)|^2 $  by (\ref{eq:var}). With the help of a change of variable, it can be shown that 
$$ I_2\lesssim \sum_{n=n_0}^{n_1-1} \left|k_h(t^{n+1}_h)-k_h(t^n_h)\right| \left|l(t^{n}_h)-l(t^{n_1}_h)\right| \lesssim [|k_h|(t)-|k_h|(s)] \sup_{s\leq t_1\leq t_2 \leq t} |l(t_2)-l(t_1)|.$$
We also have estimated the first term $I$, it remains to deal with the second one $II$ (due to (\ref{no11})). We recall that
$$ k_h(t^{n+1}_h) - k_h(t^n_h) \in \NN(C(t_h^{n+1}),\tilde{x}_h(t^{n+1}_h)).$$
The hypomonotonicity property of the proximal normal cone (see Definition \ref{def} and Proposition \ref{prop:maxhypo}) yields
\begin{align*}
-\left\langle k_h(t^{n+1}_h)-k_h(t^{n}_h), \tilde{x}_h(t^{n+1}_h)-\tilde{x}_h(t^{n_0}_h)\right\rangle  \leq & \frac{1}{2\eta} |\tilde{x}_h(t^{n+1}_h)-\tilde{x}_h(t^{n_0}_h)|^2|k_h(t^{n+1}_h)-k_h(t^{n}_h)| \\
 & + \frac{1}{2\eta} d_H(Q(t^{n_0}_h),Q(t^{n+1}_h))^2|k_h(t^{n+1}_h)-k_h(t^{n}_h)| \\
 & + \left|k_h(t^{n+1}_h)-k_h(t^{n}_h)\right| d_H(Q(t^{n_0}_h),Q(t^{n+1}_h)).
\end{align*}
Indeed, the point $\tilde{x}_h(t^{n_0}_h)\in Q(t^{n_0}_h)$ and possibly does not belong to $Q(t^{n+1}_h)$ so we have to replace it by its projection onto $Q(t^{n+1}_h)$, which makes appear the two last quantities.
We produce a similar reasoning for $\tilde{x}_h(t^{n}_h)$, in noting that 
$$ |\tilde{x}_h(t^{n}_h) - \tilde{x}_h(t^{n+1}_h)|\leq \left[v(t_h^{n+1})-v(t_h^n)\right] + |l(t_h^{n+1})-l(t_h^n)|.$$
We deduce that the second term is bounded by
\begin{align*}
 II \leq & \frac{1}{\eta} \sum_{n=n_0}^{n_1-1}  |\tilde{x}_h(t^{n+1}_h)-\tilde{x}_h(t^{n_0}_h)|^2 |k_h(t^{n+1}_h)-k_h(t^{n}_h)| e^{-[|k_h|(t^{n+1}_h)-|k_h|(t^{n_0}_h)]/\eta} \\
& + c \left[|k_h|(t)-|k_h|(s)\right]\left( \left[|v|(t)-|v|(s)\right] + \left[|v|(t)-|v|(s)\right]^2  + \sup_{s\leq t_1\leq t_2 \leq t} |l(t_2)-l(t_1)| \right),
\end{align*}
where $c$ is a numerical constant and $v$ is given by Definition \ref{def:abs} ($v$ controls the variation of the set $C(\cdot)$).\\
It remains to study the term $S_2$ corresponding to (\ref{no2}). Using
$$ \eta \frac{e^{-[|k_h|(t^{n}_h)-|k_h|(t^{n_0}_h)]/\eta}-e^{- [|k_h|(t^{n+1}_h)-|k_h|(t^{n_0}_h)]/\eta}}{|k_h|(t^{n+1}_h)-|k_h|(t^{n}_h)} \geq e^{- [|k_h|(t^{n+1}_h)-|k_h|(t^{n_0}_h)]/\eta},$$
it comes
\begin{align}  
\lefteqn{II+S_2} & & \nonumber \\
 &  &  \lesssim \left[|k_h|(t)-|k_h|(s)\right]\left( \left[|v|(t)-|v|(s)\right] + \left[|v|(t)-|v|(s)\right]^2  + \sup_{s\leq t_1\leq t_2 \leq t} |l(t_2)-l(t_1)| \right). \label{eq:imp} 
\end{align}
Finally with (\ref{eq:var}), we obtain that
\begin{align}
\lefteqn{|\tilde{x}_h(t)-\tilde{x}_h(s)|^2} & & \nonumber \\
 & & \lesssim \sup_{s\leq t_1\leq t_2 \leq t} |l(t_2)-l(t_1)|^2 + \sup_{s\leq t_1\leq t_2 \leq t} |l(t_2)-l(t_1)| + \left[|v|(t)-|v|(s)\right]+\left[|v|(t)-|v|(s)\right]^2. \label{eq:sup} 
\end{align}
So we get
\begin{align*} 
\sup_{T_m\leq s\leq t \leq T_{m+1}} |\tilde{x}_h(t)-\tilde{x}_h(s)|^2  \lesssim & \sup_{s\leq t_1\leq t_2 \leq t} |l(t_2)-l(t_1)|^2 + \sup_{s\leq t_1\leq t_2 \leq t} |l(t_2)-l(t_1)| \\
 & + \left[|v|(T_{m+1})-|v|(T_{m})\right] + \left[|v|(T_{m+1})-|v|(T_{m})\right]^2.
\end{align*}
As in Proposition \ref{prop3}, the uniform continuity of $v$ implies that $\left[|v|(T_{m+1})-|v|(T_{m})\right]$ can be assumed small with respect to $r$ (it suffices to take $\tau$ sufficiently small). Hence, the collection of indices $m$ is finite and so (\ref{eq:var}) becomes
\be{eq:var2} |k_h|(\tau) - |k_h|(0) \leq K' \ee
for some numerical constant $K'$ depending only on $l$ (by its uniform mudulus continuity and its $L^\infty$-norm).
Dividing the time-interval $\I$ with subintervals of length $\tau$, we deduce that the total variation  $|k_h|(T)$ is (uniformly with respect to $h$) bounded by a constant depending on $l\in C^0(\I)$. \\
Then (\ref{eq:var2} implies that (\ref{eq:sup}) holds for every $s,t\in\I$, which gives the uniform continuity of $\tilde{x}_h$.
\findem

\begin{rem} \label{rem:extraterm} We have detailed the bounds of the terms $I$, $II$ and $S_2$ in order to make explicit the extra term in (\ref{eq:imp}). In the continuous versions of these results with time-independent set $C(t)=C$, it is well-known that we get $ II + S_2 \leq 0$ (see (7) in Lemma 1.1 of \cite{LS} for example). \\
This new term makes appear the ``variation'' of the set $C(\cdot)$ in the corresponding interval $[s,t]$.
\end{rem}

\mb Then we deduce the following result:

\begin{thm} \label{thm:fin} Consider an admissible set-valued map $C$, varying in an absolutely continuous way and $f,\sigma$ be functions satisfying the assumptions of Theorem \ref{thm:principal2}. Let $X$ be the unique process, solution of
$$ \left\{ \begin{array}{l} 
                  dX_t + \NN(C(t),X_t) \ni f(t,X_t)dt + \sigma(t,X_t) dB_t \vsp \\
                  X_0=u_0 ,
                 \end{array} \right.$$
where $u_0\in C(0)$ is a non-stochastic initial data, given by Theorem \ref{thm:principal2}. 
For all $\omega\in \Omega$, we construct the following discretized process:
\begin{itemize}
\item for all $t\in[0,t^{1}_h]$
 $$ \tilde{X}^h_t(\omega) := u_0 ,\quad \textrm{ for all } t\in[0,t^{1}_h]$$
\item for all $t\in]t^{n}_h, t^{n+1}_h]$
 $$ \tilde{X}^h_t(\omega) : = \PPP_{C(t^{n+1}_h)}\left[\tilde{X}^h_{t^n_h}(\omega) + hf(t^n_h,\tilde{X}^h_{t^n_h}(\omega)) + \sigma(t^n_h,\tilde{X}^h_{t^n_h}(\omega)) \left(B_t-B_{t^n_h}\right) \right]. $$
\end{itemize}
Then for almost every $\omega\in \Omega$, we have:
$$ \lim_{h\rightarrow 0} \left\| \tilde{X}^h(\omega) - X(\omega)\right\|_{L^\infty(\I)} = 0.$$
\end{thm}

\mb We refer the reader to Theorem 5.1 of \cite{Saisho} (mainly Lemmas 5.1 and 5.2 of this work) for a detailed proof of such result, using compactness arguments.

\section{Example of applications with particular set-valued maps} \label{sec:particular}

In this section, we deal with some particular moving sets $C$, defined as an intersection of complements of convex sets. More precisely, we check that the assumptions (made in the previous theorems) are satisfied in the framework of \cite{JuJu}. Let us recall it.

\mb We consider the Euclidean space $\R^d$, equipped with its euclidean metric $|\ |$, its inner product $\langle\cdot,\cdot \rangle$ and $\B$ the closed unit ball in $\R^d$. Let $\I:=[0,T]$ be a bounded closed time-interval and for $i\in\{1,..,p\}$ let $g_i:\I\times \R^d \rightarrow \R$ be functions (which can be thought as ``constraints''). We introduce the sets $Q_i(t)$ for every $t\in \I$ by:
$$ Q_i(t):=\left\{x\in \R^d,\ g_i(t,x)\geq 0\right\},$$
and the following one
$$ Q(t):= \bigcap_{i=1}^p Q_i(t),$$
which represents the set of ``feasible configurations $x$''. We suppose that for all $t\in \I$ and all $i\in\{1,..,p\}$, $g_i(t,\cdot)$ is a convex function. We suppose there exists $c  >0$ and for all $ t$ in $[0,T]$ open sets $U_i(t) \supset Q_i(t) $ verifying
\begin{equation} \label{Ui}
 \tag{A0}
d_H(Q_i(t), \R^d \setminus U_i(t)) > c,
\end{equation}
where $d_H $ denotes the Hausdorff distance. \\
Moreover we assume that there exist $\alpha,\beta,M,\kappa>0$ such that $g_i\in C^2\left(\I\times (Q_i+\kappa\B) \right)$ and satisfies:
\begin{equation} \label{gradg}
\forall t\in \I, \ x\in U_i(t), \qquad \alpha \leq |\nabla_x g_i(t,x)| \leq \beta, \tag{$A1$}
\end{equation}
\begin{equation} \label{dtg}
\forall t\in \I, \ x\in U_i(t), \qquad |\partial_t g_i(t,x)| \leq \beta, \tag{$A2$}
\end{equation}
and
\begin{equation} \label{hessg}
\forall t\in \I, \ x\in U_i(t), \qquad |D^2_x g_i(t,x)| \leq M. \tag{$A3$}
\end{equation}
\begin{equation} \label{dtgradg}
\forall t\in \I, \ x\in U_i(t), \qquad |\partial_t \nabla_x g_i(t,x)| \leq M. \tag{$A4$}
\end{equation}
We denote by
$$ I(t,x):=\left\{i,\ g_i(t,x)=0\right\} $$ the set of ``active contraints'' and 
$$ I_\rho(t,x):=\left\{i,\ g_i(t,x)\leq \rho\right\}, $$
for some $\rho>0$.
As explained before, we want to deal with admissible sets $Q(t)$ and so we have to make the following important assumption: 
there exist constants $\rho,\gamma>0$ such that for all $t\in \I$ and all $x\in Q(t)$
\be{reverse}
\sum_{i\in I_\rho(t,x)} \lambda_i |\nabla_x g_i(t,x)| \leq \gamma \left| \sum_{i\in I_\rho(t,x)} \lambda_i \nabla_x g_i(t,x)\right |, \tag{$R_\rho$}
\ee
for every nonnegative coefficients $\lambda_i$.
We refer to \cite{JuJu} for a first use of this kind of ``reverse triangle inequality'' ($R_0$) and ($R_\rho$).

\mb We want to apply the previous results to the set-valued map $Q(\cdot)$ in order to get well-posedness results for the following stochastic sweeping process
\be{eqstobis} \left\{ \begin{array}{l} 
                  dX_t + \NN(Q(t),X_t) \ni f(t,X_t)dt + \sigma(t,X_t) dB_t \vsp \\
                  X_0=u_0 \in Q(0).
                 \end{array} \right. \ee


\mb By Proposition 2.8 in \cite{JuJu}, we know that for all $t\in\I$, $Q(t)$ is uniformly prox-regular and we can describe its proximal normal cone.

\begin{prop} \label{prop:prox} Under the assumption ($R_0$), there exists a constant $\eta$ such that for all $t\in \I$, the set $Q(t)$ is $\eta$-prox-regular. \\
Moreover, for all $t\in \I$ and $x\in Q(t)$,
$$ \NN(Q(t),x)=\sum_{i\in I(t,x)} \NN(Q_i(t),x) = -\sum_{i\in I(t,x)} \R^+ \nabla_x\, g_i(t,x).$$
\end{prop}


\mb In order to study the set-valued map $Q(\cdot)$, we need this technical lemma

\begin{lem} \label{gooddir}
 There exist constants $ \nu,\tau,\rho',r >0$ such that for all $t \in \I $ and $ x \in Q(t)$: there exists $u$ satisfying:
\begin{itemize}
 \item $|u|=1$
 \item for all $s\in[t-\tau,t+\tau]\cap \I$, $y\in B(x,2r)$ and  $i\in I_{\rho'}(s,y)$,
\be{eq:gooddir} \left<\nabla_x\, g_i(s,y),u \right>\geq\nu. \ee
\end{itemize}
\end{lem}

\begin{rem} Such result was already proved in Lemma 2.10 \cite{JuJu} with $r=\tau=0$.
\end{rem}

\dem Let $t\in\I$ and $x\in Q(t)$, we set the following cone
$$ \NN_\rho(Q(t),x):= -\sum_{i\in I_\rho(t,x)} \R^+ \nabla_x \, g_i(t,x)$$
and its polar cone
$$ C_\rho(Q(t),x):= \NN_\rho(Q(t),x)^\circ:=\left\{ w \in \R^d \virg \forall v \in \NN_\rho(Q(t),x) \virg \langle  v, w\rangle \leq 0 \right\}.$$
According to the classical orthogonal decomposition of a Hilbert space as the sum of mutually polar cones (see~\cite{Moreau}), we have:
$$ Id= \PPP_{\NN_\rho(Q(t),x)} + \PPP_{C_\rho(Q(t),x)},$$
where $\PPP$ denotes the Euclidean projection. \\
So let us consider for $i\in I_\rho(t,x)$ the corresponding decomposition of $\nabla_x \, g_i(t,x)$:
$$ \nabla_x \, g_i(t,x) = a_i+b_i \in \NN_\rho(Q(t),x) + C_\rho(Q(t),x). $$
Assumption (\ref{gradg}) gives us: $|b_i|\leq |\nabla_x \, g_i(t,x)| \leq \beta$.
Since $a_i \in \NN_\rho(Q(t),x) $, it can be written: $a_i = - \sum \lambda_j  \nabla_\qqq \, g_j(t,x) \virg \lambda_j \geq 0$ involving $$|b_i|= |\nabla_x \, g_i(t,x)-a_i|= \left| \sum_{j \neq i}\lambda_j  \nabla_x \, g_j(t,x)+(1+ \lambda_i) \nabla_x \, g_i(t,x) \right|.$$
Then using the inverse triangle inequality (\ref{reverse}) and Assumption (\ref{gradg}), we get:
$$ |b_i| \geq \frac{\alpha}{\gamma} \left(\sum \lambda_j +1 \right) \geq \frac{\alpha}{\gamma}.$$
As a consequence, it comes: 
\begin{equation}
 \frac{\alpha}{\gamma} \leq |b_i| \leq \beta.
\label{bi}
\end{equation}
Since $0\in \NN_\rho(Q(t),x)$ and $a_i=\PPP_{\NN_\rho(Q(t),x)}(\nabla_x \, g_i(t,x))$, we obtain
\begin{align} 2\langle b_i,-\nabla_x \, g_i(t,x) \rangle & =  |b_i-\nabla_x \, g_i(t,x)|^2-|b_i|^2 - |\nabla_x \, g_i(t,x)|^2 \nonumber \\
& =|a_i|^2 -|b_i|^2 - |\nabla_x \, g_i(t,x)|^2 \nonumber \\ 
& \leq -|b_i|^2 \leq -\frac{\alpha^2}{\gamma^2}. 
\label{eq:gooddir2} 
\end{align}
Now we set 
\begin{equation}
u:=\frac{\sum_{i\in I_\rho(t,x)} b_i}{|\sum_{i\in I_\rho(t,x)} b_i|} \in C_\rho(Q(t),x).
\label{defu}
\end{equation}
This is well-defined because (\ref{eq:gooddir2}) and Assumption (\ref{gradg}) imply that for any $j\in I_\rho(t,x)$
\begin{equation}\left|\sum_{i\in I_\rho(t,x)} b_i\right| \geq \frac{1}{\beta}\left\langle \sum_{i\in I_\rho(t,x)} b_i,\nabla_x \, g_j(t,x) \right\rangle \geq \frac{1}{\beta}\left\langle b_j, \nabla_x \, g_j(t,x)      \right\rangle            \geq \frac{\alpha^2}{2\beta \gamma^2}.
\label{sombi}
\end{equation}
Then (\ref{eq:gooddir}) for $(s,y)=(t,x)$ follows from (\ref{bi}) and (\ref{sombi}) with $$ \dsp \nu'=\frac{\alpha^2}{2 \gamma^2 p\beta}.$$
For $\rho'=\rho/2$, Lipschitz regularity (\ref{gradg}) and (\ref{dtg}) imply for $y\in B(x,2r)$ and $s\in[t-\tau,t+\tau]\cap \I$ with $r,\tau \leq \rho/(8\beta)$
$$ I_{\rho'}(s,y) \subset I_{\rho}(t,x).$$
Thanks to Assumptions (\ref{hessg}) and (\ref{dtgradg}), it can be shown
$$  \left<\nabla_x\, g_i(s,y),v \right> \geq \nu' - M(2r+\tau) \geq \nu,$$
with $\nu=\nu'/2$ and $r,\tau \leq \nu/(2M)$. \\
Consequently (\ref{eq:gooddir}) holds for $\rho'=\rho/2$, $\nu=\frac{\alpha^2}{4 \gamma^2 p\beta}$ and $r,\tau \leq \min\{ \rho/(8\beta), \nu/(2M)\}$.
\findem

\begin{prop} \label{prop1} For all $t\in \I=[0,T]$, the set-valued map $Q(\cdot)$ is admissible. 
\end{prop}

\dem The uniform prox-regularity of $Q(t)$ is asserted in Proposition \ref{prop:prox}. It also remains to check Property (\ref{hyp5}). \\
Let $t\in \I$ and $x\in Q(t)$ and $r,\tau,u$ given by the previous Lemma. We fix $s\in [t-\tau,t+\tau] \cap \I$ and $y\in B(x,2r) \cap \partial Q(s)$. Every $v\in\NN(Q(s),y)\cap {\mathbb B}$ can be written as follows
$$ v = -\sum_{i\in I(s,y)} \lambda_i \nabla_x\, g_i(s,y)$$
with nonnegative coefficients $\lambda_i$. By Assumption (\ref{gradg}),
$$\beta \sum_{i\in I(s,y)} \lambda_i \geq 1.$$ 
Then from Assumptions (\ref{hessg}) and (\ref{dtgradg}), we deduce a Lipschitz regularity for the gradient $\nabla_x\, g_i$ and so, we get
\begin{align}
\langle v,-u\rangle & = \sum_{i\in I(s,y)} \lambda_i \langle \nabla_x g_i(s,y) , u \rangle  \nonumber \\
 &  \geq \nu \sum_{i\in I(s,y)} \lambda_i \nonumber \\
 &  \geq  \nu/\beta. \label{delta}
\end{align}
We set $\delta:=\nu/\beta$ and $u_x:=-u$. \\
So we have proved that for every $x\in\partial Q(t)$ we can find $u_x$, satisfying for all $y\in \partial Q(s) \cap B(x,2r)$ with $|t-s|\leq \tau$ and $v\in \NN(Q(s),y)\cap {\mathbb B}$
\be{eq} \langle v, u_x \rangle \geq \delta. \ee
Then Property (\ref{hyp5}) is obtained by choosing a bounded covering of $\dsp \bigcup_{s,\ |s-t|\leq \tau} \partial Q(s)\subset \R^d$ with balls of radius $r$.
\findem

\begin{prop} \label{Lipschitz} The set-valued map $Q$ is Lipschitz continuous with respect to the Hausdorff distance.
\end{prop}

\mb This result was already proved in \cite{JuJu}. We give a proof for an easy reference.

\dem
Consider $t,s \in [0,T]$ and $x\in Q(t)$, let us construct a point close to $x$ belonging to $Q(s)$.
Let $u $ given by Lemma \ref{gooddir}, we introduce 
 $z(h):=x+hu$ with $h >0$. We claim that for $h< h_l:=\min(c,\rho/(\beta +\nu)) $ , 
$$ \dsp  \forall i\in \{1,...,p\}, \qquad g_i(t,z(h))\geq h \nu$$ 
(where $c $ and $\nu $ are introduced in (\ref{Ui}) and (\ref{eq:gooddir})).
Indeed for  $h< c$, due to the convexity of $ g_i(t,\cdot)$, it comes
$$
g_i(t,x+hu) \geq g_i(t,x) + h \langle \nabla_x \, g_i(t,x), \ u \rangle.
$$
 As a consequence, for $i\in I_\rho(t,x)$,
$$ \dsp g_i(t,x+hu) \geq h\nu ,$$
by (\ref{eq:gooddir}).
Moreover for every $i \notin  I_\rho(t,xq) $, according to Assumption (\ref{gradg}), $$g_i(t,xq+hu) \geq \rho -h\beta \geq h\nu $$  if $\dsp h<\frac{\rho}{\beta + \nu} $.
Thus for $h<h_l$, we have $g_i(t,x+hu)\geq h \nu$ for all $i \in \{1,..,p\}$. 
That is why we deduce from Assumption (\ref{dtg}) that $z(h)\in Q(s)$ if $ h \nu \geq \beta |t-s|$. Setting $\ell:=\dsp  \frac{\beta  h_l}{\nu} $, if $|t-s|<\ell $, it can be written that
$$ d_{Q(s)}(x) \leq \inf_{h \nu  \geq  \beta|t-s|} |x -z(h)| \leq \frac{ \beta}{\nu} |t-s|.$$
Consequently we obtain if $|t-s|<\ell $,
$$ d_H(Q(t),Q(s)) = \max\left(\sup_{x\in Q(t)} d_{Q(s)}(x), \sup_{x\in Q(s)} d_{Q(t)}(x) \right)\leq  \frac{\beta}{\nu} |t-s|.$$
This inequality is actually satisfied for any $t, s \in [0,T]$. To check it, it suffices to divide the corresponding interval into subintervals of length $\ell$ and to apply the triangle inequality. \findem

\begin{rem} We would like to show the importance of Assumption ($R_0$) for the Lipschitz regularity of $Q$ with an example. Consider the two-dimensional space $\R^2$ ($d=2$) and chose two constraints
$$ g_1(t,x):=x_2 \qquad \textrm{and} \qquad g_2(t,x):=-e^{-x_1}-x_2-t.$$
\begin{figure} 
\centering
\begin{tabular}{cc}
\psfrag{x1}[l]{$x_1$}
\psfrag{x2}[l]{$x_2$}
\psfrag{Q}[l]{$Q(0)$}
\includegraphics[width=0.4\textwidth]{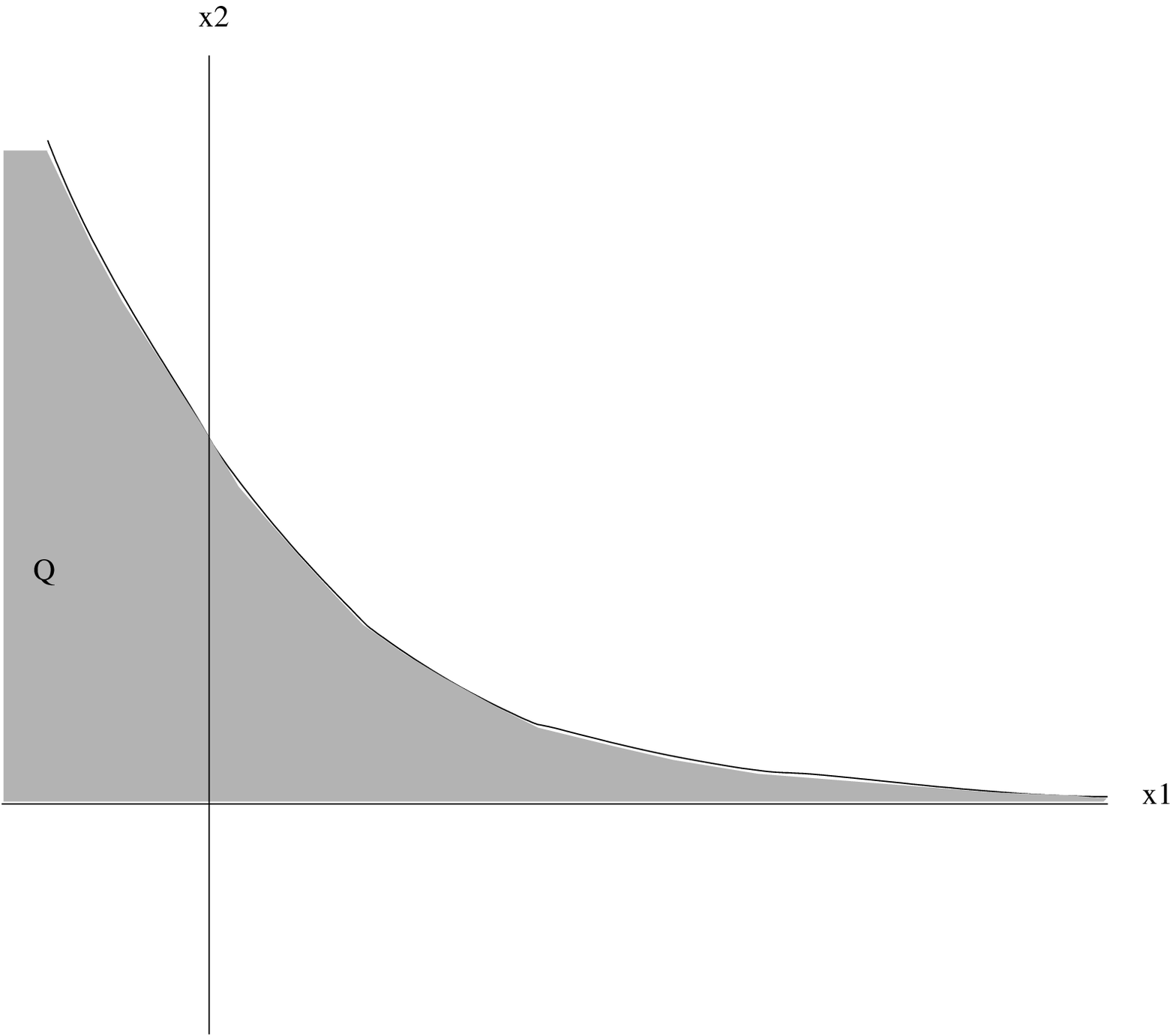} 
 \quad & 
\psfrag{x1}[l]{$x_1$}
\psfrag{x2}[l]{$x_2$}
\psfrag{Q}[l]{$Q(t)$}
\psfrag{t}[l]{$-t$}
\psfrag{c}[l]{$-log(t)$}
\includegraphics[width=0.4\textwidth]{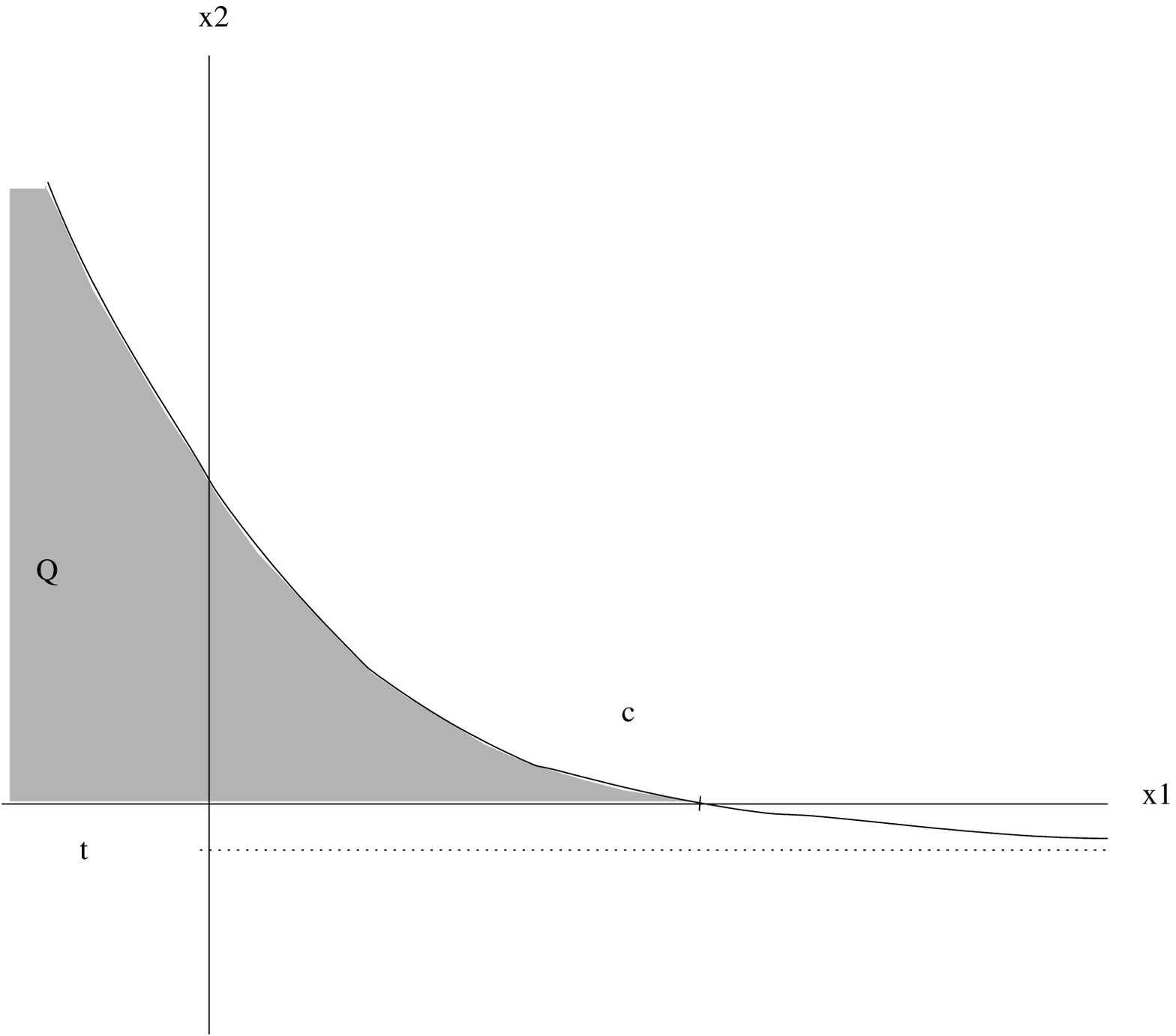} 
\end{tabular}
\caption{The sets $Q(0)$ and $Q(t)$ for some $t>0$.}
\label{fig:Q}
\end{figure}
\noindent The two constraints are smooth and convex and satisfy the assumptions (\ref{gradg})-(\ref{dtgradg}).
Then it is easy to check that for all $t\geq 0$ the line $\{x_2=0\}$ belongs to the set $Q(t)$ (see Figure \ref{fig:Q}). However for $t>0$, we have
$$ Q(t) \subset \left\{(x_1,x_2),\ x_2\geq 0 \textrm{ and } x_1\leq -\log(t) \right\}.$$
So we deduce that $d_H(Q(t),Q(0)) = \infty$ for all $t>0$.
Computing the gradients, we get
$$ \nabla_x\, g_1(t,x) = (0,1) \qquad \textrm{ and } \qquad \nabla_x\, g_2(t,x) = (-e^{-x_1},-1).$$
Assumption ($R_0$) is not satisfied as we have:
$$ \left|\nabla_x\, g_1(t,x) + \nabla_x\, g_2(t,x) \right| = e^{-x_1}$$
which could be as small as we want.
\end{rem}

\gb{ \bf Conclusion:} In this framework (Assumptions (\ref{Ui})-(\ref{dtgradg}) and (\ref{reverse})), we have proved that the set-valued map $Q$ is admissible and Lipschitz continuous.\\
We can also apply Theorem \ref{thm:principal2} and obtain well-posedness results concerning the following stochastic differential inclusion (with maps $f$ and $\sigma$ which are bounded and Lipschitz with respect to the second variable):
$$ 
\left\{
\begin{array}{l}
 \dsp dX_t + \NN(Q(t),X_t) \ni f(t,X_t)dt + \sigma(t,X_t)dB_t \vsp \\
 X_t\in Q(t) \vsp \\
 X_0=u_0 ,
\end{array}
\right. 
$$
for all $u_0\in Q(0)$. Moreover, an Euler scheme can be used to approach the solution, as it is convergent thanks to Theorem \ref{thm:fin}.

\gb This numerical scheme (\ref{schema}) makes appear the Euclidean projection $\PPP_{Q(t)}$ onto the set $Q(t)$. As the set $Q(t)$ is not assumed to be convex but only uniformly prox-regular, we have not efficient numerical algorithm in order to compute this projection. \\
In the case of the set $Q$, defined (see the beginning of this section) as the intersection of complements of convex sets, an implementable scheme has been proposed in \cite{Maurygrain, esaim, MV} and then extended to moving sets $Q(t)$ in \cite{JuJu}. The idea is to replace $Q(t)$ with a convex set $\tilde{Q}(t,x)$ (depending on the variable $x$) defined for any point $x\in Q(t)$ as follows:
$$
\Qc(t,x)  := \left\{ y\in \R^d \virg g_{i}(t,x) + \langle \nabla g_{i}(t,x), 
y - x \rangle \geq 0 \quad \forall \, i \right\}.
$$
This set is convex and included into $Q(t)$, due to the convexity of the functions $g_i(t,\cdot)$.
This substitution is convenient because classical methods can be employed to compute the projection onto a convex set. Yet this replacement raises some difficulties for the numerical analysis which are solved in the framework of deterministic sweeping process, in proving that $\Qc(t,x)$ is a good local approximation of $Q(t)$ around the point $x$ (we refer the reader to \cite{JuJu} for a detailed study and corresponding results).
The new algorithm (instead of (\ref{schema})) is the following one: for $l\in C^0(\I)$
we now define a discretized solution $\tilde{x}_h$ as follows:
\be{schema2} \left\{ \begin{array}{l}
            \tilde{x}_h(t) := u_0 ,\quad \textrm{ for all } t\in[0,t^{1}_h] \vsp \\
            \tilde{x}_h(t) : = \PPP_{\Qc(t^{n+1}_h,\tilde{x}_h(t^n_h))}\left[\tilde{x}_h(t^n_h) + l(t) - l(t^n_h) \right] ,\quad \textrm{ for all } t\in]t^{n}_h, t^{n+1}_h].
          \end{array}
\right.  
\ee
It would be interesting to prove its convergence.
In \cite{JuJu}, the second author has already proved its convergence in the framework of sweeping process for $l\in W^{1,1}(\I)$. The extension of this result to continuous functions $l\in C^0(\I)$ is still open. Indeed, the main difficulty rests on the lack of regularity of the map $x\rightarrow \Qc(t,x)$ (even for time-independent constraints $g_i$).

\gb We finish this work by briefly presenting a stochastic model of crowd motion, which is an extension of a deterministic one introduced in \cite{MV2}. We refer the reader to~\cite{TheseJu, MV, MV2} for a complete and detailed description of this model, which takes into account the direct conflict between people. 

\mb We consider $N$ persons identified to rigid disks. The
center of the $i$-th disk is denoted by $\qq_i\in \R^2$ and its radius by $r_i$. Since overlapping is forbidden, the  
vector of positions $\qqq=(\qq_1,..,\qq_N) \in
\R^{2N}$ has to belong to the ``set of feasible configurations'', defined by 
\be{eq:Q} Q:=\left\{ \qqq  \in \R^{2N},\ D_{ij}( \qqq) \geq 0\quad   \forall \,i \neq j  \right\},  \ee 
where $ D_{ij}(\qqq)=|\qq_i-\qq_j|-(r_i+r_j) $ is the signed distance between disks $i$ and $j$. 

\mb By denoting by $\UU(\qqq)=(U_1(\qq_1),..,U_N(\qq_N)) \in \R^{2N}$ the global spontaneous velocity of the individuals, the crowd motion model can be written
$$ d\qqq + \NN(Q,\qqq) \ni \UU(\qqq)dt.$$

\mb Now people's hesitation or panic can be modelled by a stochastic perturbation $\sigma(t,\qqq)dB_t$, which gives the following stochastic model:
\be{eq:model3} d\qqq + \NN(Q,\qqq) \ni \UU(\qqq)dt + \sigma(t,\qqq)dB_t,\ee
associated to a Brownian motion $(B_t)_{t>0}$. \\
In order to get well-posedness results, it also suffices to check that $Q$ is admissible. This was already proved in \cite{Saisho2} by a direct approach and later in \cite{JuJu} via showing Assumption ($R_\rho$) (since Assumptions (\ref{Ui}), (\ref{gradg}) and (\ref{hessg}) can be easily satisfied).


\begin{thm} If $\UU,\sigma(t,\cdot)$ are Lipschitzean and $\sigma$ is bounded, then (\ref{eq:model3}) is well-posed: for every initial condition $\qqq(0)\in Q$, there is a unique solution to (\ref{eq:model3}).
\end{thm}

\begin{rem} Since we deal with a constant set $Q$, \cite{Saisho} already allows to conclude. However the current results allow to get the well-posedness of (\ref{eq:model3}) when the radii depend on time under the following assumptions: $r_i$ are uniformly time-Lipschitz and
$$ \inf_{t\in[0,T]} \inf_{i} r_i(t) >0.$$
\end{rem}

\bibliographystyle{plain}
\bibliography{biblio}

\end{document}